\documentclass[11pt]{article}
\usepackage{amsmath}
\usepackage{amsthm}
\usepackage{verbatim}
\usepackage{color}
\usepackage{subfigure}
\usepackage{graphicx}
\usepackage{booktabs}
\def\RR{I\!\!R}
\oddsidemargin=.5cm \evensidemargin=0mm \newfont{\titre}{cmbx12 at 16pt}

\newtheorem{thm}{Theorem}[section]

\newtheorem{rem}{Remark}[section]
\small

\def\mathr{{\sf I\kern-0.1em R}}
\def\Mathz{{\Large \kern0.3em\sf \Large Z\kern-0.4em Z\kern0.3em}}
\def\mathz{{\kern0.3em\sf Z\kern-0.8em Z\kern0.2em}}
\def\mathc{{\kern0.1em\sf C\kern-0.8em C\kern0.2em}}

\newcommand*{\ang}[1]{\left\langle #1 \right\rangle}

\def\11{I\!\!I}  \def\RR{I\!\!R} 
\def\RR{{\sf I\kern-0.1em R}}

\newcounter{mnote}
\title{Finite determining parameters feedback control for distributed nonlinear dissipative systems - a computational study \\
}
\author{{Evelyn ~Lunasin$^1$ and Edriss S. ~Titi} $^{2, 3}$}
\date{\small  $^{1}$ Department of Mathematics, United States Naval Academy, Annapolis, MD 21402\\
{\em lunasin@usna.edu}\\
$^{2}$ Departments of Mathematics, Texas A\&M University, College Station,
TX 77843-3368, USA.\\
$^{3}$  Department of Computer Science and Applied Mathematics,
The Weizmann Institute of Science, Rehovot 76100, Israel.\\
{\em titi@math.tamu.edu} and {\em edriss.titi@weizmann.ac.il}
\\
July 8, 2017}

\begin{document}
\maketitle{}

\abstract{

We investigate the effectiveness of a simple finite-dimensional feedback control scheme for globally stabilizing solutions of infinite-dimensional dissipative evolution equations introduced by Azouani and Titi in \cite{AzTi13}.   This feedback control algorithm overcomes some of the major difficulties in control of multi-scale processes: It does not require the presence of separation of scales nor does it assume the existence of a finite-dimensional globally invariant inertial manifold.  In this work we present a theoretical framework for a control algorithm which allows us to give a systematic stability analysis,  and present the parameter regime where stabilization or  control objective is attained.  In addition, the number of observables and controllers that were derived analytically and implemented in our numerical studies  is consistent with  the finite number of determining modes that are relevant to the underlying physical system.  We verify the results computationally in the context of the Chafee-Infante reaction-diffusion equation, the Kuramoto-Sivashinsky equation, and other applied control problems, and observe that the control strategy is robust and independent of the model equation describing the dissipative system.

  {\bf Keywords.} {\em Globally stabilizing feedback control, Chafee-Infante, Kuramoto-Sivashinsky, reaction-diffusion, Navier-Stokes equations, feedback control, data assimilation, determining modes, determining
 nodes, determining volume elements.}\\
{\bf Mathematics Subject Classification (2000):} 35K57, 37L25, 37L30, 37N35, 93B52, 93C20, 93D15.

\tableofcontents


\section{Introduction}

The main characteristic of infinite-dimensional dissipative  dynamical systems is that they posses finite number of degrees of freedom. There is a long history of how to take advantage of this fact in order to design reduced finite-dimensional feedback control for stabilizing their solutions (see, e.g., \cite{Armou_Christofides2000, Cao-Kevrekidis-Titi-2001, Christofides_Armou2000, Yannis-Control} and references therein). In \cite{Yannis-Control}, for instance, the authors survey and employ various reduction methods, such as the proper orthogonal decomposition and the nonlinear Galerkin method (approximate inertial manifolds tools), to design feedback control of various reaction-diffusion systems.  Moreover, other techniques for finite-dimensional control and  global stabilization of the Kuramoto-Sivashinsky equations  (KSE) can be found in \cite{Armou_Christofides2000, Christofides_Armou2000}.

In \cite{Armou_Christofides2000}, the authors employed nonlinear Galerkin method to design a nonlinear low dimensional output feedback controllers for the KSE and as well as for the Korteweg-de-Vries-Burgers equations.  The reduced order model is used for control synthesis using techniques from geometric control methods. Moreover,  rigorous numerical criteria for stabilizing solutions of  the Navier-Stokes equations by linear feedback control, through the Galerkin and nonlinear Galerkin methods, were investigated in \cite{Cao-Kevrekidis-Titi-2001}.

Similarly, the authors in \cite{Christofides_Armou2000}  derived a simple linear feedback used to achieve  stabilization of the zero solution for any value for the instability parameter, for any initial conditions.   Their feedback control algorithm starts by deriving an ODE approximation of the integrated form of the controlled KSE system, subject to periodic boundary conditions, that captures  the dynamics of the unstable modes.   Necessary and sufficient conditions under which  a linear static output feedback controller that uses five control actuators and $2s+1$ measurements (where $2s$ is the number of unstable modes in the system) were  derived for global exponential stabilization of the zero solution of the KSE.

In \cite{DMEC_2006} the authors addressed numerous issues  pertaining to  predictive control of linear parabolic PDEs with constraints in states and control actuators.  Their method also starts with the construction of reduced order models that captures the dominant dynamics of the infinite-dimensional system and then they apply the state constraints in the reduced order model.  Moreover, they also compared the performance of the predictive control strategy for several different modal projections.

Efficient algorithms for control and stabilization of infinite-dimensional dynamical systems is a topic of significant interest in  several areas of science and engineering.   Some major areas of application include stabilizing flame front propagation, chemical process control  (see, e.g.,  \cite{Christofidesbook} for application of control for  catalytic rod, coating processes and nonlinear control of Czochralski crystal growth processes), control in biological systems and  control of turbulence/chaos in lasers \cite{ILG2007}.      For most of these applications, the model dynamics are known to exhibit low dimensional spatio-temporal chaos due to their dissipative behavior.   Thus, active monitoring and control of the state of the system takes advantage of the reduced order models (cf. \cite{Yannis-Control} and references therein).  Additionally, in several of these works the construction of the control design are based on the local  linearization of the reduced equations about a reference steady state solution.     Examples of valuable results can be found in \cite{Ahuja2009, Christofides2000,Christofides2003,Foias-Jolly-Kevrekidis-Sell-Titi,Jolly-Kevrekidis-Titi, Li_Christofides_2007, Ros:rd, Ros:rd-im, Yannis-Control} and references therein.  Despite continuing efforts among many researchers  there are still many analytical issues left unresolved.  In particular, the real-time implementation of  the control algorithm in industrial control systems is observed to be a nontrivial issue.

In this paper we aim to present a unified rigorous approach, with complete proofs and without any {\it ad hoc} uncheckable assumptions, for stabilizing solutions of the KSE and other dissipative systems finitely many parameters (see also \cite{Kalantarov_Titi2015, Kalantarov_Titi2017} for similar approach concerning the nonlinear damped wave equations and other related systems).
In particular, we investigate a global feedback control  algorithm that side-tracks some of the difficulties mentioned earlier.   Our current study stems from a recent analytical result for a simple finite-dimensional feedback control algorithm  introduced in  \cite{AzTi13}  for globally stabilizing solutions of infinite-dimensional dissipative evolution equations.  We note that the feedback controller proposed in \cite{AzTi13} can be used for stabilization of the zero steady state solutions of a wide class of nonlinear dissipative wave equations as well (see, e.g. \cite{Kalantarov_Titi2015} and references therein). The newly proposed global feedback control scheme uses finitely many observables and controllers which is consistent and stems from the fact that such systems possess a finite number of determining parameters or degrees of freedom, for example, finite number of determining Fourier modes, determining nodes, and determining interpolants and projections.  As it has been noted in \cite{AzTi13}, this proposed algorithm does not require the presence of separation in spatial scales, in particular, it does not assume the existence of a globally invariant inertial manifold (see, e.g. \cite{Ros:rd-im} for using inertial manifold for feedback control).   See also \cite{Foias-Jolly-Kravchenko-Titi1, Foias-Jolly-Kravchenko-Titi2, JST, JST2} for relevant work on  {\it determining forms}.

Furthermore, the subsequent implementation and analysis of  this simple control algorithm as ``nudging'' form of data assimilation, for example, in \cite{Azouani-Olson-Titi, FJT14, FLT14, FLT_horizontal, FLT_PG, FLT_porous, MTT2015}, has resulted in analytical progress that has potential valuable impact in numerous fields.  These works present rigorous analytical support to continuous in time data assimilation algorithm for certain dynamical systems that can identify the full state of the system knowing only coarse spatial mesh observational measurements, not of full state variables of the model, but only of some selected state variables in the system.  Recently,  an improved continuous data assimilation  algorithm for the 3D Planetary Geostrophic model that requires observations of the temperature only is presented in \cite{FLT_PG} (see also \cite{Ghil1977, Ghil-Halem-Atlas1978}). This provides a rigorous justification to an earlier conjecture of Charney  in \cite{Charney1969} which states that temperature history of the atmosphere, for certain simple atmospheric models, will determine other state variables.  Earlier implementations of the ``nudging" algorithm suffers from a  lack of concrete theoretical analysis which left physicists searching for the optimal or suitable nudging coefficient (or relaxation parameter) through expensive numerical experiments.   The works mentioned above provided some stepping stones to rigorous justification to some of the earlier conjectures of experts in the field.    In addition, the systematic theoretical framework of the proposed global control scheme allowed the authors in these works to provide sufficient conditions on the spatial resolution of the collected spatial coarse mesh data and the relaxation parameter that guarantees that the approximating solution obtained from this algorithm converges to the unknown reference solution over time (with the assumption that the observational data measurements are free of noise).  The computational study implementing these algorithms under drastically more relaxed conditions can be found in  \cite{Altaf} and \cite{Gesho} .

 It is also worthwhile to note that the authors in \cite{Bessaih-Olson-Titi} studied the performance of this  feedback control algorithm applied to data assimilations when the observational data contains stochastic measurement errors.   The algorithm is applied to the 2D Navier-Stokes equations paradigm. The resulting equation in the algorithm is a Navier-Stokes-like equation with stochastic feedback term that attunes the large scales in the approximate solution to those of the  reference solution corresponding to the measurements.  The authors in \cite{Bessaih-Olson-Titi}  found resolution conditions on the observational data and the relaxation/nudging parameter in which the expected value of the difference between the approximate solution, recovered by proposed linear feedback control algorithm and the exact solution is bounded by a factor which depends on the Grashof (Reynolds) number multiplied by the variance of the noise, asymptotically in time.  Moreover, an extension of this algorithm for discrete spatiotemporal data, and error analyses has recently been reported in \cite{FMTi}.   In addition, statistical solutions (i.e., probability invariant measures) are also established in \cite{FMTi} based on discrete data. More recently in \cite{M_T}, the authors obtain uniform in time estimates for the error between the numerical approximation given by the Post-Processing Galerkin method of the downscaling algorithm and the reference solution, for the 2D NSE.   Notably, this uniform in time error estimates provide a strong evidence for the practical reliable numerical implementation of this algorithm.

In this article we verify the validity of this global feedback control algorithm using various numerical experiments.   In particular, we implement the control scheme for a simple reaction-diffusion equation, the Chafee-Infante equation, which is the real version of the complex  Ginzburg-Landau equation, and for the  Kuramoto-Sivashinksy equation (KSE). (See, e.g.,  \cite{Gomes_2016A, Gomes_2016B, Thompson_2016} for a few recently proposed control strategy for stabilizing falling liquid film flows in the stochastic and deterministic case).  We give rigorous global stability analysis for the feedback control algorithm for stabilizing a reference unstable steady state solution for the one-dimensional KSE and derive the sufficient conditions on the parameter regime to achieve the control objective.  We also implement  the  algorithm for control of nonlinear parabolic system governing model equations for catalytic rod with or without uncertainty, which are examples of case studies studied in  \cite{Christofidesbook}.    We show that the simpler feedback control algorithm performs as expected achieving similar success.



\section{Various types of interpolant operators for feedback controllers}\label{s:interpolant}
This preliminary section serves two purposes, one is to  enumerate the set of interpolant operators used in the stability analysis presented
\cite{AzTi13},  and second is to set the notation and serve as an introduction for our  computational studies and analysis.  For $\varphi \in H^1([0, L])$ we define
\begin{equation}
\|\varphi\|_{H^1}^2:=\frac{1}{L^2} \, \int_{0}^{L}
\varphi^2(x)\, dx + \int_{0}^{L} \varphi^2_{x} (x) \, dx.
\label{norm-h1}
\end{equation}
Consider a general linear map $I_h: H^1([0,L])\rightarrow L^2([0,L])$ which is an interpolant operator that approximates the identity operator with error of order $h$.  Specifically, it approximates the inclusion map $i: H^1\hookrightarrow  L^2$ such that the estimate
\begin{equation}\label{IH-operator}
\|\varphi - I_h(\varphi)\|_{L^2}\leq c h \|\varphi\|_{H^1},
\end{equation}
holds for every $\varphi\in H^1([0,L])$, where $c$ is a dimensionless constant independent of $\phi$ and $h$.  Examples of approximate interpolant operators, discussed in \cite{AzTi13}, that we will consider here, with the general mapping property (see e.g., \cite{Jones92B, Jones93}) mentioned above are as follows:
\begin{enumerate}
\item {\it Fourier modes.}  Consider a periodic function $\varphi\in L^1[0,L]$.  The interpolant operator $I_h$ acting on  $\varphi\in L^1[0,L]$ is defined as  the projection onto the first  $N$ Fourier modes;
\begin{equation}\label{ih-modal}
I_h(\varphi)=\frac{a_0}{2} + \sum_{k=1}^{N} {a_k }\, \cos{\frac{k\pi x}{L}} + \sum_{k=1}^{N} {b_k }\, \sin{\frac{k\pi x}{L}}, \quad h=\frac{L}{N},
\end{equation}
where the Fourier coefficients are given by $$a_k=\frac{2}{L} \int_{0}^{L} \, \varphi(x) \, \cos{\frac{k\pi x}{L}} \,dx,  \quad b_k=\frac{2}{L} \int_{0}^{L} \, \phi(x) \, \sin{\frac{k\pi x}{L}} \,dx. $$

\item {\it Finite volume elements.}  The volume element operator is also an interpolant operator satisfying \eqref{IH-operator}.  Given $\varphi\in L^1[0,L]$, we define
\begin{eqnarray}
I_h(\varphi) =
\sum_{k=1}^{N} \overline{\varphi}_k \; \chi_{J_k\strut} (x), \label{ih-fv}
\end{eqnarray}
where $J_k=\left[(k-1) \frac{L}{N}, k \frac{L}{N}\right),$
for $k=1,\dots,N-1$, and $J_N = \left[(N-1) \frac{L}{N},  L \right]$, $\chi_{J_k \strut}(x)$
is the characteristic function of the interval $J_k$, for $k=1,\dots,N$, serving as the actuator shape function,  and where
$$
\bar{\varphi}_k=\frac{1}{|J_k|} \, \int_{J_k} \, \varphi(x) \, \, dx=
\frac{N}{L} \, \int_{J_k} \, \varphi(x) \, \, dx,
$$
represents the amplitude for the given actuator.
Here, the local averages, $\overline{\varphi}_k$, for
$k=1,..., N,$ are the observables, which also serve as the feedback
controllers.
\item {\it Interpolant operator based on nodal values}.  Let $\varphi\in H^1[0,L]$.   We also consider here the situation where the observables are the values $\varphi(x_k^*)$, where $x_k^*\in J_k = [(k-1)\frac{L}{N}, k\frac{L}{N}]$, $k = 1, 2, \dots, N$.
In this case the interpolant operator is given by

\begin{eqnarray}
I_h(\varphi) = \sum_{k=1}^N \, \varphi(x_k^*) \,\chi_{J_k\strut} (x), \quad x\in[0, L]. \label{ih_nodal-fv}
\end{eqnarray}
where again $\chi_{J_k \strut}(x)$
denotes the characteristic function of the interval $J_k$, for $k=1,\dots,N$, serving as the actuator shape function.

\end{enumerate}

\section{The Chafee-Infante equation}

We also recall the motivating model equation used to demonstrate in its simplest case the proposed algorithm in \cite{AzTi13};  the Chafee-Infante reaction-diffusion equation (or the real version of the complex Ginzburg-Landau equation) on the interval $[0, L]$,  with no flux boundary condition is given by

\begin{subequations}
\begin{eqnarray} \label{RD}
\frac{\partial u}{ \partial t}-\nu \, u_{xx}-\alpha u+u^3=0\\
\label{bc1}
u_x (0)=u_x (L)=0,
\end{eqnarray}
\end{subequations}
for $\alpha > 0$, large enough, and for a given initial condition $u(x,0) = u_0(x)$.


The following scheme was proposed as a feedback control system for
(\ref{RD})-(\ref{bc1})
to globally stabilize the steady state solution $\mathbf{v}\equiv 0$:
\begin{subequations}
\begin{eqnarray}
\frac{\partial u}{\partial t}- \nu \, u_{xx}-\alpha u+u^3 =-\mu
I_h(u) \label{av-RD} \\
u_x(0)=u_x(L)=0, \label{av-bc1}
\end{eqnarray}
\end{subequations}
where $I_h$ is specified in section \ref{s:interpolant}. The results concerning global existence, uniqueness and stabilization for general family of finite-dimensional feedback control system that includes system (\ref{av-RD})-(\ref{av-bc1}) as
a particular case, were established in \cite{AzTi13}.   They  showed that every solution
$u$ of (\ref{av-RD})-(\ref{av-bc1}) tends to zero at an exponential rate, as
$t \rightarrow \infty$, under specific explicit estimates on
$N$, in terms of the physical parameters $\nu, \alpha, L$ and $  \mu$.  Their main global stability results are stated in the next theorem. 

\begin{thm}
\label{thm:fin-vol}
Let $N$ and $\mu$  be large enough such that $N > \sqrt{\frac{L^2\alpha}{4\pi^2\nu}}$ and $\mu  \geq \nu\left(\frac{2 \pi}{h}\right)^2> \alpha$,
where $\alpha >0$ and $h=\frac{L}{N}$. Then $\|u(t)\|_{L^2}$
tends to zero, as $t \rightarrow \infty$, for every solution $u(t)$
of (\ref{av-RD})-(\ref{av-bc1}).
\end{thm}

\subsection*{Remark 2.1}
\label{B3}
We also recall an important observation noted in \cite{AzTi13}. The assumptions used to establish Theorem \ref{thm:fin-vol} in \cite{AzTi13}, in particular, that
$N > \sqrt{\frac{L^2 \, \alpha}{4 \pi^2 \nu}}$,  has a subsequent physical justification.   The assumption is known to be consistent with the previously established dimension of the unstable manifold about
the steady state solution $\mathbf{v} \equiv 0$ of \eqref{RD}-\eqref{bc1} which is of order of $\sqrt{\frac{L^2 \, \alpha}{\nu}}$ (see for e.g., \cite{BabVish92, Hale88, Tem88}). In addition, there is nothing special about the zero solution. And that one can use the same idea to globally stabilize any other given solution, $v(x,t)$, of \eqref{RD}-\eqref{bc1} by using a slightly modified  feedback control in the right-hand side of \eqref{av-RD}-\eqref{av-bc1} of the form $-\mu
\sum_{k=1}^{N} (\overline{u}_k -\overline{v}_k)\; \chi_{J_k\strut} (x)$.

One important open problem  mentioned in \cite{AzTi13} is whether one can design a feedback control scheme that stabilizes the zero stationary solution using only two controllers at these nodes.  This open problem is motivated by result established in \cite{Kukavica} (see also generalization in \cite{Collet}, \cite{OliverTiti_2001}) that the dissipative system \eqref{RD}-\eqref{bc1} has two determining close enough nodes.  Our current results and efforts leave this issue as an open problem.

\subsection{Numerical results for the Chafee-Infante equation}


The Chafee-Infante equation \eqref{RD}-\eqref{bc1} has an unstable trivial steady state, for $\alpha > \frac{4\pi^2\nu}{L^2}$. Starting with some initial condition close to zero, the system encounters large scale instabilities causing the solution to grow in time until it reaches another steady solution.
  Basic numerical simulation of the 1D Chafee-Infante equation demonstrates these instabilities.   We establish below that the proposed global feedback control strategy described in \eqref{av-RD}-\eqref{av-bc1} prevents this instability.


We implement numerically the proposed  control algorithm \eqref{av-RD}-\eqref{av-bc1} with $I_h$ defined as in \eqref{ih-fv} to test a simple case.   We assume that measurements of the state $u(x,t)$ are available at the discretized positions and discrete times.   The actuator and sensor locations are distributed uniformly throughout the domain.  To give a clear illustration, in Figure 1,  we only show the results for $t\in$ $[0, 0.1]$.  The number of controls (denoted by $NC$) that we have used in the numerically simulation is  $NC = 10$.  This  is consistent with the number of unstable modes $\sqrt{L^2\alpha/\nu} = 10$ with the given parameters $\alpha=100$, $\nu=1$, $L=1$.  The value $\mu=300$ used in the simulation is much smaller than the value stated in Theorem \ref{thm:fin-vol}.  This implies that a less stringent condition is actually required than derived in the theory (which is also consistent with the implementation of the analogue data assimilation algorithm studied  computationally in \cite{Altaf} and \cite{Gesho}) .  We summarize our parameters for the two numerical test in Table 1.

\begin{figure}[h!]  
\centering																				
\parbox{5cm}{
\includegraphics[width=5.5cm]{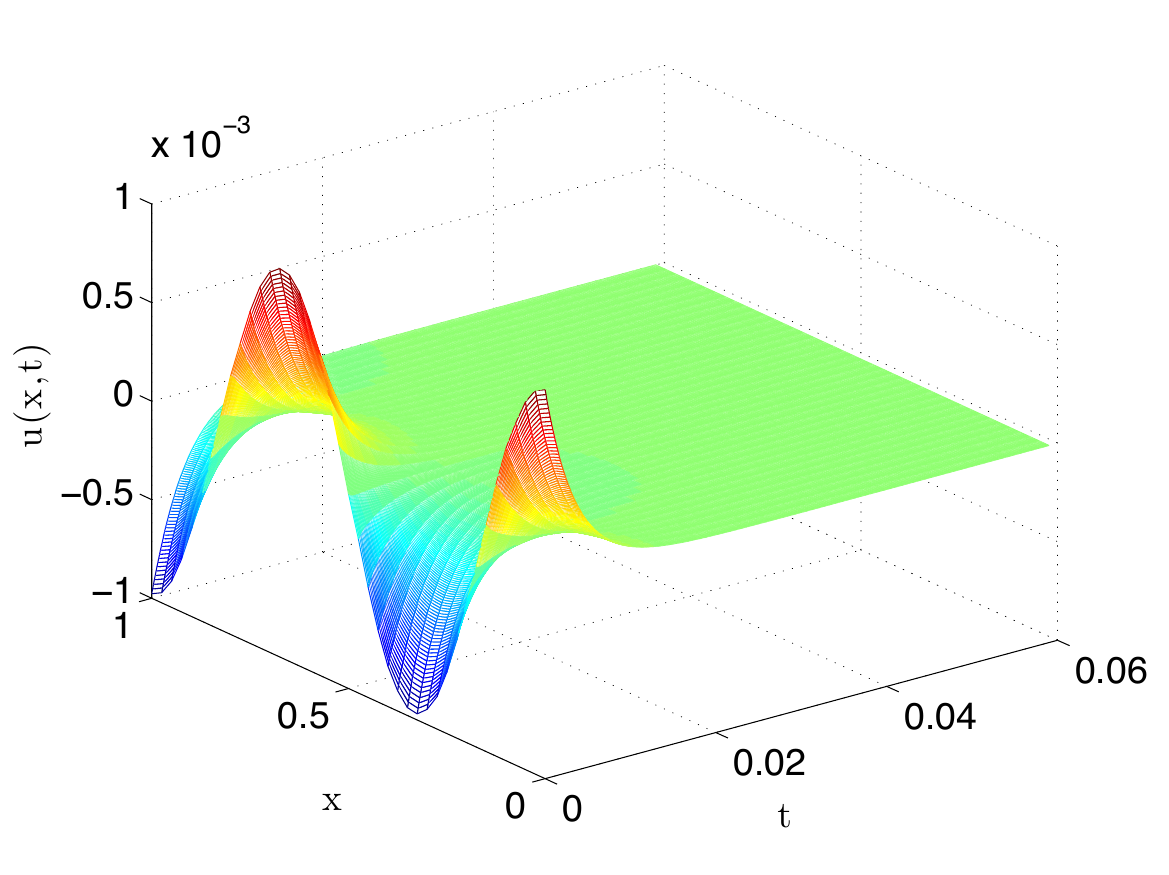}
\label{CH2A} 																			
}
\qquad
\begin{minipage}{5cm} 																	
\includegraphics[width=5.5cm]{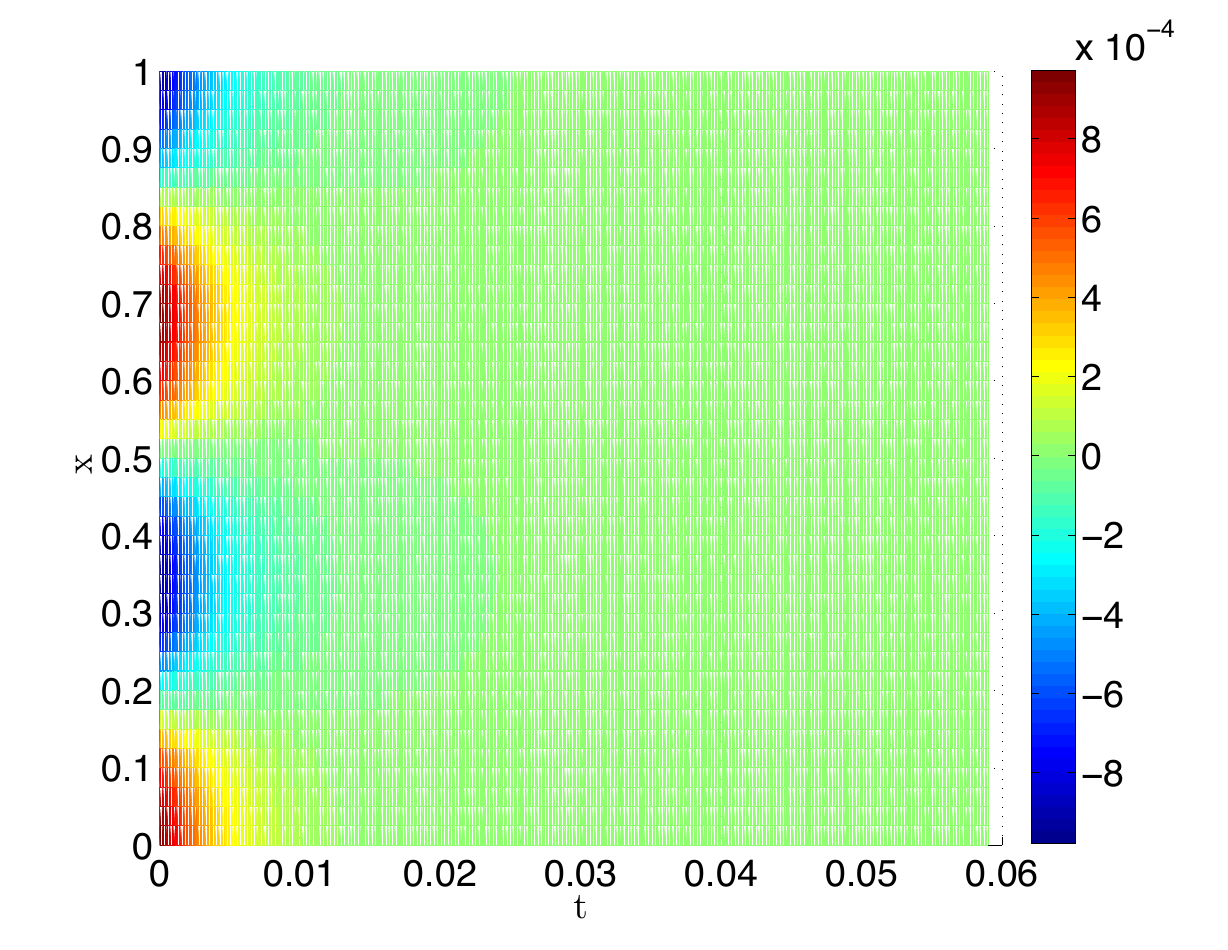}
\label{CH2B}																			
\end{minipage}
\caption{(a) Closed-loop profile showing stability of the $u(x,t)=0$ steady state solution. (b)  Top-view.}
\end{figure}	

\begin{table}
\centering
    \begin{tabular}{cclllllll}
    \toprule
    Figure & \# Actuators & $\mu$ & $\nu$ & $\alpha$ & Interpolant operator\\
    \midrule
    1 & 10 &  300 &  1 &   100 & {\it finite volume elements} \\
    \bottomrule
    \end{tabular}
    \caption{Model parameters and type of interpolant operator for the controlled and uncontrolled 1D Chafee-Infante equations}
    \end{table}

We show in the next section the performance of the proposed algorithm for other types of well-studied control problems.
We verify our results in the context the Kuramoto-Sivashinky equation, and previously studied applied control problems and observe that the proposed control strategy in \cite{AzTi13} is robust and independent of the model equation describing the dissipative system.


\section{The Kuramoto-Sivashinsky equation }

 We illustrate the application of the proposed feedback control algorithm to the global  stabilization of the zero solution of another nonlinear dissipative PDE known as the  1D KSE.   
Introduced as a turbulence model for flame front propagation by Sivashinsky in 1977 (\cite{Sivashinksy1977}), KSE also models the motion of a thin viscous fluid or thin film flowing down an inclined wall derived by Chen and Chang in 1986 (\cite{Chen_Chang_1986}).    This equation is also known to model magnetized plasmas as derived in \cite{Cohen1976, LaQuey1975} and as a model for  chemical reaction-diffusion processes as derived in \cite{Kuramoto-Tsusuki1974}.
  In one space dimension,  it is written as

\begin{eqnarray}\label{KSE}
\frac{\partial u}{\partial t}= - \gamma\dfrac{\partial^2u}{\partial x^2} - \nu \, \dfrac{\partial^4u}{\partial x^4}- u\dfrac{\partial u}{\partial x}, \quad x\in[0,L], \label{KSE}
\end{eqnarray}
\noindent
subject to the periodic boundary conditions,
and initial condition:
\begin{equation}\label{KSE-init}
u(x,0) = u_0(x),
\end{equation}
where $u(x,t)$,  for example, describes the height of the film fluctuations, and the parameters $\gamma$ and $\nu$ are given positive constants.  We assume that $\int_0^L u_0(x) dx=0$ which implies that $\int_0^L u(x,t)dx = 0$.   Equation \eqref{KSE} can be nondimensionalized by substituting $u \rightarrow \gamma u /\tilde{L}, \; t\rightarrow t\tilde{L}^2/\gamma, \; x\rightarrow \tilde{L}x,$ and $\nu \rightarrow \tilde{L}^2 \gamma\nu$, with $\tilde{L} = \dfrac{L}{2\pi}$.  In this case one gets the same equation as before with the modification $\gamma = 1$ and $L=2\pi$.

To be more precise, after rescaling the equation becomes
\begin{eqnarray}\label{E1}
\frac{\partial u}{\partial t}= - \dfrac{\partial^2u}{\partial x^2} - \nu \, \dfrac{\partial^4u}{\partial x^4}- u\dfrac{\partial u}{\partial x}, \quad x\in[0,2\pi].
\end{eqnarray}

Rigorous analytical studies have revealed that the KSE enjoys finite-dimensional asymptotic (in time) behavior  (see, e.g.,  \cite{CJT,Cock97, Cos:nse, Cos88, FMRT, Foias-Sell-Titi, Hale88, NLT_85, Robinson, Sell-You, Tem88},
 and references therein). Dissipative systems, such as the KSE, possess  finite-dimensional global attractors (\cite{BabVish92, Cos:nse, Cos88, Robinson,Sell-You,Tem88}), and finite number  of determining modes (\cite{Foias-Prodi, Foias-Manley-Temam-Treve, FMRT, Jones93}), determining nodes (\cite{FMRT, Foias-Temam-nodes1, Foias-Temam-nodes2, Foias91, Jo:nse2D, Jones93, Kukavica}), determining volume
elements (\cite{Foias91,Jones92B}) and other finite number of determining parameters (degrees of freedom)  such as finite elements and other interpolation polynomials (\cite{CJT,Cock97,Foias-Temam-nodes2}.)  The KSE also enjoys the property of separation of spatial scales,  which guarantees the existence of a
finite-dimensional globally invariant inertial manifolds (see, e.g., \cite{Cos:nse, Cos88, Foias88, Foias-Sell-Titi, Tem88}, and references therein), although it is not needed for the implementation of the  control algorithm we propose here.  The effectiveness of the feedback control strategy proposed in \cite{AzTi13} relies on the fact that the instabilities in such systems occurs solely at large spatial scales, and hence all that needed is to control these large spatial scales.


To give some justification for the choice of the number of  feedback controllers (for our method the actuators and sensors are in the same locations) it is necessary to know the number of unstable modes of the steady state ${\bf v}(x)\equiv 0$ for a given bifurcation parameter value $\nu$.  In order to obtain a heuristic value for the number of unstable modes,  one can perform a simple analysis by linearizing
equation \eqref{E1} about ${\bf v} \equiv 0$, subject to periodic boundary conditions, to obtain the linear equation

\begin{eqnarray}
\frac{\partial \bf v}{\partial t}= - \dfrac{\partial^2 \bf v}{\partial x^2}- \nu \, \dfrac{\partial^4\bf v}{\partial x^4} , \quad x\in[0,2\pi] \label{KSE-lin}.
\end{eqnarray}
\noindent
Assuming a particular solution of the form  $\mathbf{v}(x, t)=a_k(t) \, e^{ikx}$  yields the equation
\begin{equation}\label{e:coeff}
\dot{a_k}= (k^2 - \nu k^4) a_k,
\end{equation}
for the time dependent coefficients.  We solve equation \eqref{e:coeff} with initial conditions $a_0(x)=A_k \in \RR$,    $k \in \mathbf{Z}\setminus\left\{0\right\}$, which yields
\begin{equation*}
a_k(t)=A_k \, e^{k^2(1 - \nu k^2) t}.
\end{equation*}
This shows that all the low wave numbers $k <\frac{1}{\sqrt{\nu}}$ are unstable.  Thus,
one needs at least $\frac{2}{\sqrt{\nu}}$ number of parameters  to stabilize  $\mathbf{v} \equiv 0$ and that  the nonlinear system \eqref{KSE},  is locally unstable when $\nu<1$.

\subsection{The KSE without feedback control}


 To show some  proof-of-concept corresponding to our heuristic calculations of the number of unstable modes, we present our simulations for the case where the $\nu > 1$   and the case where the $\nu < 1$.   We choose an initial condition $1e^{-10} * \cos x \left(1 + \sin x\right)$ for both cases.  Observe that for $\nu=1.1 > 1$, our linear stability analysis shows exponential decay to the linearly stable steady state zero solution. Figure 2a is consistent with this result.    The final profile of the film height is given in Figure 2b.

\begin{figure}[h!]  \label{KSE0A2}	
\centering																				
\parbox{5cm}{
\includegraphics[width=5.5cm]{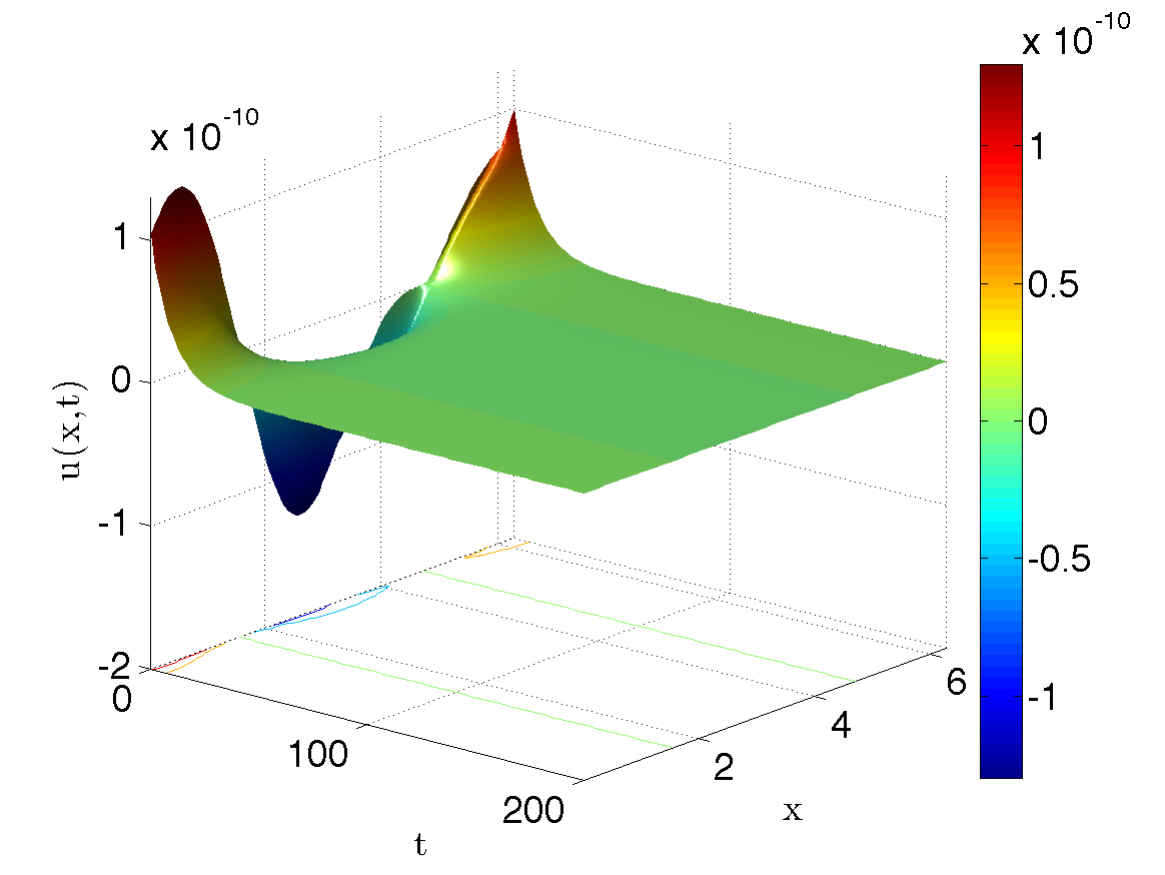}
}
\qquad
\begin{minipage}{5cm} 																	
\includegraphics[width=5.5cm]{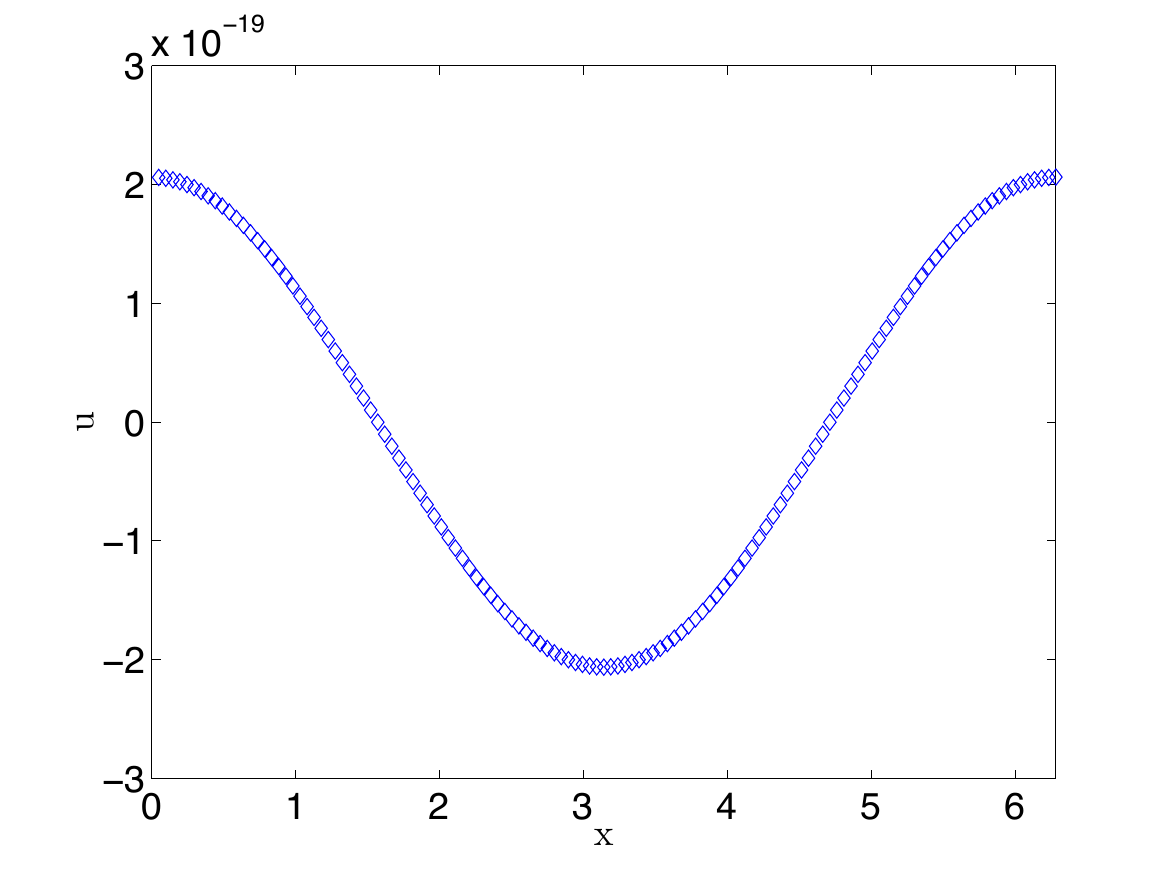}
\end{minipage}
\caption{ (a) Open-loop profile showing stability of the $u(x,t)=0$ steady state solution when $\nu = 1.1 > 1$ (b)  Profile of $u(x, t=200)$. }

\end{figure}	

For the case where $\nu = 4/15 < 1$, our numerical simulation is illustrated in Figures 3a and 3b.  In the context of a thin film flowing in an inclined surface, Figures 3a and 3b  illustrate the unwanted wavy fluctuations that develop in time.  Small perturbation of the film height at $t=0$ results in unwanted  structures starting to form at $t = 32$ and then more precisely the solution move towards a stable traveling wave pattern starting around $t=80$.    We show only up to some particular final time to display a clear transition between patterns or structures.

\begin{figure}[h!]  
\centering																				
\parbox{5cm}{
\includegraphics[width=5.5cm]{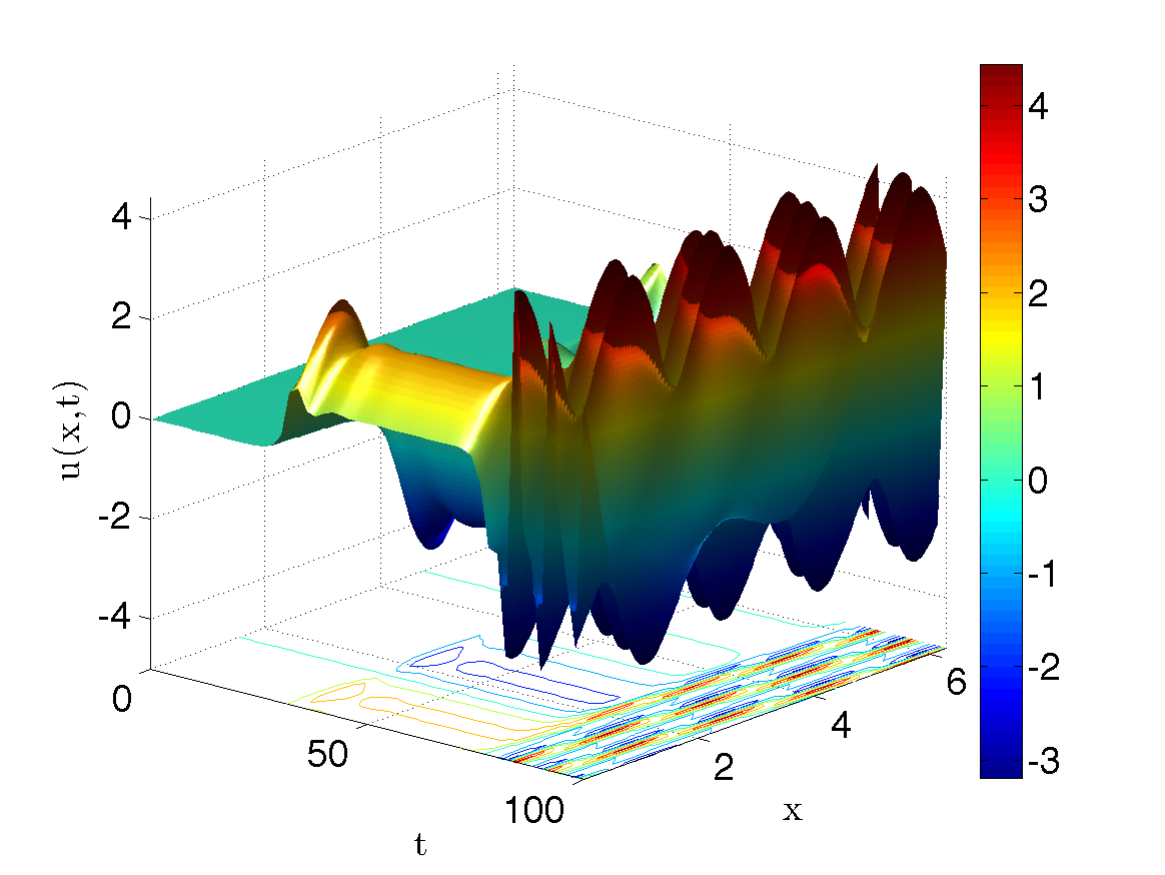}
\label{KSEOB1} 																			
}
\qquad
\begin{minipage}{5cm} 																	
\includegraphics[width=5.5cm]{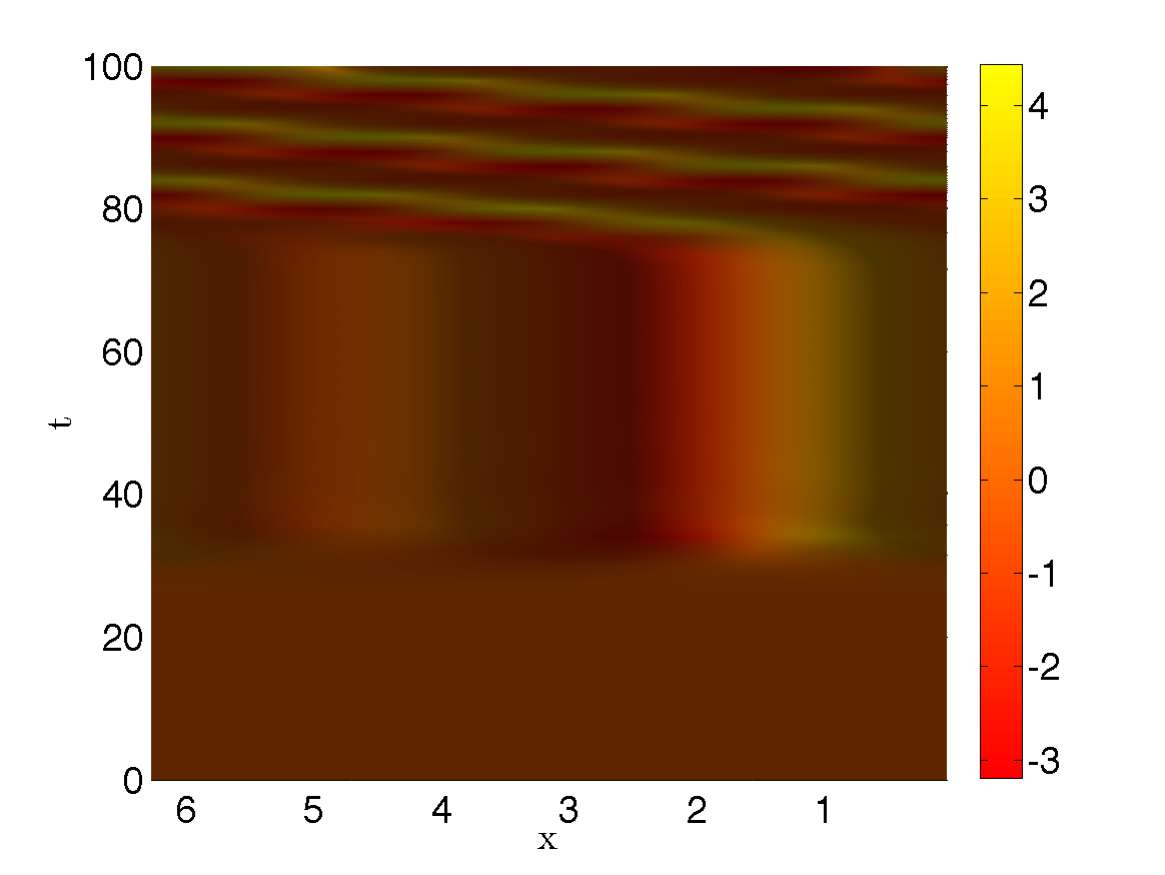}
\label{KSE0B2}																			
\end{minipage}
\caption{(a) Open-loop profile showing instability of the $u(x,t)=0$ steady state solution when $\nu = 4/15  < 1$. (b)  Top view profile of $u(x, t)$.}
\end{figure}	

A common goal (see, e.g., \cite{Ahuja2009}, \cite{Christofides2000}, \cite{Lee_Tran_2005}) is to identify a control strategy to maintain a uniform thin film's thickness by automating an actuator such as a blower or suction that serves as a source or sink at specific locations in the inclined plane where the liquid film is flowing.  One also needs to design feedback control laws using a broad class of actuators that are likely to achieve the stabilization of the film height to the zero solution in real time making the implementation tractable for industrial control system problems.


 To stabilize the steady state solution ${\bf v}=0$ of \eqref{E1}, we implement the proposed feedback control algorithm in \cite{AzTi13}, applied to KSE, which is given as follows:
\begin{eqnarray}
\frac{\partial u}{\partial t}= - \dfrac{\partial^2u}{\partial x^2} - \nu \, \dfrac{\partial^4u}{\partial x^4}- u\dfrac{\partial u}{\partial x} - \mu \,I_h(u), \quad x\in[0,2\pi] \label{KSE-ih}
\end{eqnarray}
\noindent
subject to the periodic boundary conditions, and initial condition $u(x,0) = u_0(x)$, with $\int_0^{2\pi}u(x,0)\, dx = 0$, where the
interpolant operator $I_h$ acting on $u$  can be defined as one of the general interpolants listed in section \ref{s:interpolant} satisfying certain properties.    Note that the interpolant polynomial here is shifted by its average in the whole domain $[0, 2\pi]$ to maintain the invariance of the zero spatial average for the controlled equation,  (see e.g., equation \eqref{ihhh} in the numerical results section below).  We begin by establishing global well-posedness and stability  results for the proposed feedback control algorithm.

\subsection{Existence, uniqueness and global stability results for the control of  solution of the KSE (general case).} \label{sec:exis-uni}

We recall that the existence and uniqueness of solution,  as well as the existence of finite-dimensional global attractor to the system \eqref{KSE} (or equivalently \eqref{E1}) were first established in \cite{NLT_85} (see also \cite{Tem88}) for odd initial data.  The long-time boundedness results was later improved and extended to any mean-zero initial data in \cite{Collet, Goodman}.  Improved estimates were later established in \cite{Bronski2006} and \cite{GiacomeliOtto2005}.

In this section we will establish the global existence, uniqueness and stability results for the general feedback algorithm stated as

\begin{eqnarray}
\frac{\partial u}{\partial t} +  \dfrac{\partial^2u}{\partial x^2} + \nu \, \dfrac{\partial^4u}{\partial x^4}+u\dfrac{\partial u}{\partial x} =- \mu \,(I_h(u) - I_h(u^*)),\quad x\in [0, 2\pi]
 \label{geu22app}
\end{eqnarray}
where $I_h$ is a feedback control interpolant operator that is  stabilizing a specific  solution $v = u^*$ of  \eqref{E1}. This will be accomplished by choosing $\mu$ large enough and then choosing $h$ small enough, under the assumptions \eqref{e:assump2}, below, and that

\begin{equation}
\mu > \dfrac{4}{\nu}\mbox{\quad  and\quad }
\nu \geq  \mu \, ch^4.\quad
\label{assump.366app}
\end{equation}
Under the conditions stated in \eqref{assump.366app} and \eqref{e:assump2} we derive formal {\it a-priori} bounds for the difference $u-u^\ast$ that are
essential for guaranteeing the decay of $\|u - u^\ast\|_{L^2}$ to zero at an exponential rate..  These {\it a-priori} estimates, together with the global existence of $u^\ast$,  form  the key elements for showing the global existence of the solution $u$ of \eqref{geu22app}.  Notably, these formal steps can be justified by the Galerkin approximation procedure.  Uniqueness is obtained using similar energy type estimates, as shown below.

To show convergence of $u$ to the reference solution $u^\ast$, we consider the difference $w = u-u^\ast$.  The evolution equation for $w$, which is obtained by subtracting \eqref{E1} and \eqref{geu22app}, is
\begin{equation}\label{e:diff}
w_t +\nu w_{xxxx} + w_{xx}=  -ww_x -u^\ast w_x -   wu_x^\ast -\mu\, I_h(w).
\end{equation}
Multiplying by $w$, integrating by parts with respect to $x$  using periodic boundary
conditions (the cubic nonlinear term disappears) and using

\begin{equation}\label{IH-split}
I_h(w) \, w= \left(I_h(w)-w \right)w + w^2,
\end{equation}

we obtain
\begin{eqnarray}
\begin{aligned}
&\frac{1}{2} \frac{d} {d t}  \int_{0}^{2\pi} w^2 \, dx + \nu \, \int_{0}^{2\pi} w^2_{xx} \, dx\\
& = \|w_x\|^2_{L^2} + \frac{1}{2} \int_{0}^{2\pi} \, w^2 \, u_x^\ast \, dx - \mu\|w \|^2_{L^2} + \mu \,  \int_{0}^{2\pi} \left(w - I_h(w)\right) \, w \ dx.
\end{aligned}
\end{eqnarray}
A straightforward calculations, using the Poincar\'e, Cauchy-Schwarz, Young and Agmon $(\|\varphi\|^2_{L^\infty} \leq \|\varphi\|_{L^2} \|\varphi_x\|_{L^2})$ inequalities, 
 yield
\begin{eqnarray}
\begin{aligned}
&\frac{1}{2} \frac{d} {d t}  \|w \|^2_{L^2} + \nu\, \|w_{xx}\|^{2}_{L^2} \\&\leq
\, \frac{2}{\nu} \|w\|^2_{L^2}+ \frac{\nu}{8}\|w_{xx}\|^2_{L^2}+  \frac{1}{2} \|w\|^2_{L^2} \, \| u_x^\ast\|_{L^\infty} -\mu\|w\|^2_{L^2}    +\mu\, \| I_h(w)-w\|_{L^2}\|w\|_{L^2}   \\
&\leq \frac{\nu}{8}\|w_{xx}\|^2_{L^2}+\left(\frac{2}{\nu}-\mu +\| u_x^\ast\|_{L^\infty} \right)\|w\|^2_{L^2} + \mu c h  \|w_{x}\|_{L^2}  \|w\|_{L^2}\\
&\leq \frac{\nu}{8}\|w_{xx}\|^2_{L^2}+\left(\frac{2}{\nu}-\mu +\| u_x^\ast\|_{L^\infty} \right)\|w\|^2_{L^2} + \mu c h  \|w\|^{3/2}_{L^2}  \|w_{xx}\|^{1/2}_{L^2}\\
&\leq \frac{\nu}{8}\|w_{xx}\|^2_{L^2}+\frac{\nu}{4}\|w_{xx}\|^2_{L^2}+\left(\frac{2}{\nu}-\mu +\| u_x^\ast\|_{L^\infty} \right)\|w\|^2_{L^2} + \frac{3}{4}\left(\frac{\mu^4 ch^4}{\nu}\right)^{1/3}\|w\|^{2}_{L^2}\\
& =: Q.
\end{aligned}
\end{eqnarray}
Hereafter, we abuse the notation for an arbitrary constant $c$, which may  change from line to line.  By assumption (\ref{assump.366app}) that $c\mu h^4 \leq \nu$ we have  $\left(\dfrac{\mu^4 ch^4}{\nu}\right)^{1/3} \leq \dfrac{\mu}{4}$.  
Thus,

\begin{eqnarray}
\aligned
Q &\leq \frac{3\nu}{8}\|w_{xx}\|^2_{L^2}+\left(\frac{2}{\nu}- \frac{3}{4}\mu +\| u_x^\ast\|_{L^\infty} \right)\|w\|^2_{L^2} \\
 & \leq  \frac{3\nu}{8}\|w_{xx}\|^2_{L^2}+\left(-  \frac{\mu}{4} +\| u_x^\ast\|_{L^\infty} \right)\|w\|^2_{L^2},
\endaligned
\end{eqnarray}
where the last inequality is due to the  assumption  (\ref{assump.366app}), i.e.,  $\mu > 4/\nu$.
\noindent
In conclusion we have
$$
\frac{1}{2} \frac{d} {d t}  \|w \|^2_{L^2}+ \frac{5}{8}\nu \, \|w_{xx}\|^{2}_{L^2} \leq \left(-  \frac{\mu}{4} +\| u_x^\ast\|_{L^\infty} \right)\|w\|^2_{L^2},
$$
and by Gronwall's inequality it follows that
\begin{equation}
\|w(t)\|_{L^2}^2 \leq e^{\left(-\frac{\mu}{2} + \frac{2}{t}\int_0^t  \| u_x^\ast (s)\|_{L^\infty}\, ds\right) t} \, \|w(0)\|_{L^2}^2.
\label{E10}
\end{equation}
\noindent
Since $u^\ast$ is a solution of KSE, it suffices to show that for large enough $t$ we have
$$-\frac{\mu}{2} + \frac{2}{t}\int_0^t  \| u_x^\ast (\cdot, s)\|_{L^\infty}\, ds < -\frac{\mu}{8}$$
to get  an upper estimate, that decays like $\displaystyle e^{-\frac{\mu}{8}t}$, for $\|w(t)\|_{L^2}^2$ .  This can be obtained by requiring additional assumption on $\mu$ in terms of the size of the absorbing ball for KSE (\cite{NLT_85},\cite{Tem88}).  We then proceed noting that by  the one dimensional interpolation inequality $ \| u_x^\ast\|_{L^\infty} \leq  \|u_{x}^\ast \|^{1/2}_{L^2}\|u_{xx}^\ast \|^{1/2}_{L^2}$ and by Poincar\'e inequality we have $ \| u_x^\ast\|_{L^\infty} \leq \left(\frac{1}{2\pi}\right)^{1/2} \|u_{xx}^\ast \|^{1/2}_{L^2}$.   From the mathematical theory of 1D KSE (see, e.g., \cite{Bronski2006, Collet, GiacomeliOtto2005, NLT_85, Otto_2009, Tem88})   we have that
\begin{equation}\label{limsup1app}
\limsup_{t\rightarrow \infty}\frac{1}{t}\int_0^t  \| u_{xx}^\ast (\cdot, s)\|^2_{L^2}\, ds = R_2^2
\end{equation}
where $R_2 $ is a number which depends on the $\nu$ and $L$.  Therefore

\begin{eqnarray}\label{limsup2app}
\aligned
\limsup_{t\rightarrow \infty}\frac{1}{t}& \int_0^t  \| u_{x}^\ast (\cdot, s)\|_{L^\infty}\, ds \\
&\leq \left(\frac{1}{2\pi}\right)^{1/2}  \limsup_{t\rightarrow \infty}\left(\frac{1}{t}\int_0^t  \| u_{xx}^\ast (\cdot, s)\|^2_{L^2}\, ds\right)^{1/2} \leq \left(\frac{1}{2\pi}\right)^{1/2}  R_2.
\endaligned
\end{eqnarray}
Thus, if we assume that
\begin{equation}\label{e:assump2}
\dfrac{\mu}{8} \geq \left(\frac{2}{\pi}\right)^{1/2} \,R_2,
\end{equation}
 then
\begin{equation}
\limsup_{t\rightarrow\infty} \|w(t)\|^2_{L^2} \leq \limsup_{t\rightarrow \infty} e^{\left(-\frac{\mu}{4} + \frac{2}{t}\int_0^t  \| u_x^\ast (s)\|_{L^\infty}\, ds\right) t}  \|w(0)\|^2_{L^2} = 0.
\end{equation}

We obtain the following global stability result:

\begin{thm}
\label{thm:7.1App}
Let $\mu, \nu$ and $h$ be positive parameters satisfying assumption \eqref{assump.366app} and \eqref{e:assump2};
and that
$I_h$ satisfies \eqref{IH-operator}  with $\int_0^{2\pi} I_h(u)(x) \, dx=0$. Then for every $T > 0$ and $u_0 \in L^2([0, 2\pi]),$  system
\eqref{geu22app}  has a unique solution $u \in C([0, T], L^2) \cap L^2([0, T], H^2)$,
which also depends continuously on the initial data. Moreover,
$$\lim_{t \rightarrow \infty} \|u (t) - u^*\|_{L^2}^2 =0;
$$
and for every $\tau >0$
$$\lim_{t \rightarrow \infty} \int_t^{t+\tau}\|u_{xx}(s) - u^*_{xx}(s)\|_{L^2}^2 \, ds =0.$$
\end{thm}

The global stability in $H^1$ can be obtained by a slight modification of the previous  analysis.   Here we conclude  that under  the assumption \eqref{assump.366app} and \eqref{e:assump2},  the feedback control interpolant operator $I_h$ is stabilizing the solution $v\equiv u^*$ of \eqref{E1} in the $H^1$ norm and consequently in the $L^\infty$ norm, thanks to the one dimensional Sobolev imbedding theorem.  
Taking the inner product of \eqref{e:diff} with $-w_{xx}$, writing $I_h(u)$ as in \eqref{IH-split} and integrating by parts, yields
\begin{eqnarray}
\begin{aligned}
\frac{1}{2}& \frac{d} {d t}  \int_{0}^{2\pi} w_x ^2 \, dx +  \nu \,
\int_{0}^{2\pi} w_{xxx} ^2 \, dx \\
&= - \int_{0}^{2\pi} w_xw_{xxx} \, dx   -\frac{1}{2} \int_0^{2\pi} u_x^*(w_x )^2\, dx -\frac{1}{2}  \int_0^{2\pi}  w^2 w_{xxx}\, dx \\
&+\int_0^{2\pi} u_x^\ast (ww_{xx}) \; dx  + \mu \,
\int_{0}^{2\pi} (I_h(w)-w) \, w_{xx} \, dx - \mu\int_0^{2\pi} w_x^2 \, dx\\
\end{aligned}
\end{eqnarray}

\noindent
Applying Cauchy-Schwarz, H\"{o}lder's, Agmon and
 Gagliardo-Nirenberg interpolation
$(\|\phi_{xx}\|^2_{L^2}\leq \|\phi_x\|_{L^2}\|\phi_{xxx}\|_{L^2})$ inequalities,  and then applying Young inequality, we obtain
\begin{eqnarray}
\begin{aligned}
\frac{1}{2} \frac{d} {dt} & \|w_x\|^2_{L^2}+ \nu \|w_{xxx}\|^2_{L^2}\\
&\leq \|w_x\|_{L^2} \|w_{xxx}\|_{L^2}  + \frac{1}{2} \|u_x^* \|_{L^\infty}\|w_x\|^2_{L^2}   + \frac{1}{2}\| w\|_{L^\infty}\|w\|_{L^2} \|w_{xxx}\|_{L^2},\\
 &    \quad +  \|u_x^* \|_{L^\infty}\|w\|_{L^2} \|w_{xx}\|_{L^2}
 +  \mu \|I_h(w) - w\|_{L^2}\|w_{xx}\|_{L^2}-\mu\|w_x\|^2_{L^2}\\
  &= : Q_1 + Q_2 + Q_3 + Q_4 + Q_5 + Q_6,
 \end{aligned}
 \end{eqnarray}

with
 \begin{eqnarray}
 \begin{aligned}
\quad Q_1 &= \|w_x\|_{L^2} \|w_{xxx}\|_{L^2}  \leq  \dfrac{1}{\nu}\|w_x\|^2_{L^2}+\dfrac{\nu}{4}\|w_{xxx}\|^2_{L^2},\\
\quad Q_2  &= \frac{1}{2} \|u_x^* \|_{L^\infty}\|w_x\|^2_{L^2},      \\
Q_3 &= \frac{1}{2}\| w\|_{L^\infty}\|w\|_{L^2} \|w_{xxx}\|_{L^2} \leq \frac{1}{2}\| w\|^{3/2}_{L^2}\|w_x\|^{1/2}_{L^2} \|w_{xxx}\|_{L^2} \\ & \leq  \dfrac{1}{4\nu}\|w\|^{3}_{L^2}\|w_x\|_{L^2}+\dfrac{\nu}{4}\|w_{xxx}\|^2_{L^2}    \leq \frac{1}{2\mu\nu^2} \|w\|^6_{L^2} + \frac{\mu}{16}\|w_x\|^2_{L^2} + \dfrac{\nu}{4}\|w_{xxx}\|^2_{L^2},\\
Q_4 &=  \|u_x^* \|_{L^\infty}\|w\|_{L^2} \|w_{xx}\|_{L^2} \leq \|u_x^* \|^2_{L^\infty} \frac{h^2}{4\nu}  \|w\|^2_{L^2} + \frac{\nu}{h^2}\|w_{xx}\|^2_{L^2}\\
Q_5 &= \mu \|I_h(w) - w\|_{L^2}\|w_{xx}\|_{L^2} \ \leq  \mu c h\|w_x\|_{L^2}\|w_{xx}\|_{L^2}  \\
& \leq \frac{\mu}{16}\|w_x\|^2_{L^2} + \mu ch^2 \|w_{xx}\|^2_{L^2}\leq \frac{\mu}{16}\|w_x\|^2_{L^2} +  \frac{\nu}{h^2}\|w_{xx}\|^2_{L^2},\\
Q_6 & = -\mu\|w_x\|^2_{L^2},
\end{aligned}
\end{eqnarray}
\noindent
where for $Q_5$, we used the approximation identity given in \eqref{IH-operator}
and assumption (\ref{assump.366app}) that $\nu \geq  \mu \, ch^4.$   Collecting like terms we have
\begin{eqnarray}
\aligned
\frac{1}{2} \frac{d} {dt}  \|w_x\|^2_{L^2}+ \frac{\nu}{2} \|w_{xxx}\|^2_{L^2}
& \leq \left( \frac{1}{\nu}  - \mu + \frac{\mu}{8} +  \|u_x^* \|_{L^\infty}  \right) \|w_x\|^2_{L^2}\\
&+ \frac{1}{2\mu\nu^2}  \|w\|^6_{L^2} + \frac{h^2}{4\nu} \|u_x^\ast\|^2_{L^\infty} \|w\|^2_{L^2} + \frac{2\nu}{h^2}\|w_{xx}\|^2_{L^2}.
\endaligned
\end{eqnarray}

\noindent
By assumption (\ref{assump.366app}), that $\mu > \dfrac{4}{\nu}$, we have

\begin{eqnarray}
\aligned
\frac{d} {dt}  \|w_x\|^2_{L^2}+ \nu \|w_{xxx}\|^2_{L^2}  &
 \leq 2\left( - \frac{\mu}{2} + \|u_x^* \|_{L^\infty} \right) \|w_x\|^2_{L^2}\\
&+ \frac{1}{\mu\nu^2}  \|w\|^6_{L^2} + \frac{h^2}{2\nu} \|u_x^\ast\|^2_{L^\infty} \|w\|^2_{L^2} + \frac{4\nu}{h^2}\|w_{xx}\|^2_{L^2}.
\endaligned
\end{eqnarray}

Let $\epsilon >0$ be a given arbitrarily small number.     Thanks to Theorem \ref{thm:7.1App} there exists a $T_0 (\epsilon) > 0$, large enough, such that for all $t \geq T_0$ we have
\begin{equation}
\|w(t)\|_{L^2} < \epsilon, \quad \displaystyle \int_t^{t+\tau} \|w_{xx}(s)\|^2_{L^2} ds < K \epsilon, \quad \mbox{ and } \quad \|u_x^\ast\|_{L^\infty} \leq 2 \left(\frac{1}{2\pi}\right)^{1/2}R_2.
\end{equation}

Thus, for $t\geq T_0 + \tau$ we have

\begin{eqnarray}
\aligned
\frac{d} {dt}  \|w_x\|^2_{L^2}+ \nu \|w_{xxx}\|^2_{L^2}  &
 \leq 2\left( - \frac{\mu}{2} +  \|u_x^* \|_{L^\infty}  \right)\|w_x\|^2_{L^2}+ \\
&+ \frac{1}{\mu\nu^2}  \epsilon^6 + \frac{2h^2}{\nu} R_2^2\, \epsilon^2+ \frac{4\nu}{h^2}\|w_{xx}\|^2_{L^2}.
\endaligned
\end{eqnarray}

By assumption \eqref{e:assump2} and   Gronwall inequality we have

\begin{eqnarray}
\aligned
\|w_x(t)\|_{L^2}^2 &\leq e^{-\frac{\mu}{2}(t-T_0) } \, \|w_x(T_0)\|_{L^2}^2 + \frac{2}{\mu}\left(\frac{1}{\mu\nu^2}  \epsilon^6 + \frac{2h^2}{\nu} R_2^2\, \epsilon^2\right)\left(e^{-\frac{\mu}{2}T_0} - e^{-\frac{\mu}{2}t}\right)\\
& + \frac{4\nu}{h^2} \int_{T_0}^t e^{-\frac{\mu}{2}(t-s)}\|w_{xx}(s)\|^2_{L^2}\, ds.
\endaligned
\label{App1011}
\end{eqnarray}

Let us treat the integral $\displaystyle  J = \int_{T_0}^t e^{-\frac{\mu}{2} (t-s)   } \|w_{xx}(s)\|^2_{L^2}\, ds $.      There exists a natural number $N$ such that
$T_0 + N\tau \leq t < T_0 + (N+1)\tau$.  Therefore,

\begin{eqnarray}
\aligned
J &\leq  \int_{T_0}^{T_0 + (N+1)\tau} e^{-\frac{\mu}{2} (t-s)   } \|w_{xx}(s)\|^2_{L^2}\, ds\\
& = \sum_{k = 0}^N \int_{T_0 + k\tau}^{T_0 + (k+1)\tau} e^{-\frac{\mu}{2} (t-s)   } \|w_{xx}(s)\|^2_{L^2}\, ds\\
& \leq \sum_{k = 0}^N  e^{- \frac{\mu}{2} (t - (k+1)\tau - T_0)     }    \int_{T_0 + k\tau}^{T_0 + (k+1)\tau} \|w_{xx}(s)\|^2_{L^2}\, ds\\
&\leq K\epsilon\, e^{-\frac{\mu}{2} (t- T_0)   } \sum_{k=0}^N e^{\frac{\mu}{2}(k+1)\tau}\\
&\leq K\epsilon\, e^{-\frac{\mu}{2} (t- T_0)   }\left(\frac{ e^ {\frac{\mu}{2}(N+1)\tau} - 1}{e^{\frac{\mu}{2}\tau} - 1}                \right) e^{\frac{\mu}{2}\tau}
\endaligned
\end{eqnarray}
\noindent
Notice that $\displaystyle e^{-\frac{\mu}{2} (t- T_0)  }e^ {\frac{\mu}{2}(N+1)\tau} = e^{\frac{\mu}{2}(T_0 + (N+1) \tau - t)   } \leq e^{\frac{\mu}{2}\tau} $ and
$\displaystyle e^{-\frac{\mu}{2} (t- T_0)  }   \geq 0.$
\noindent
Thus, $J \leq K\epsilon \left( \frac{e^{\mu\tau}}    {e^{\frac{\mu}{2}\tau} - 1}     \right)$.  Using the upper bound for $J$ and noting that $e^{-\frac{\mu}{2}T_0} - e^{-\frac{\mu}{2}t} \leq 1$, we obtain from \eqref{App1011}, that for $t\geq T_0 + \tau$,

\begin{eqnarray*}
\|w_x(t)\|_{L^2}^2 &\leq e^{-\frac{\mu}{2}(t-T_0) } \, \|w_x(T_0)\|_{L^2}^2 + \frac{2}{\mu}\left(\frac{1}{\mu\nu^2}  \epsilon^6 + \frac{2h^2}{\nu} R_2^2\, \epsilon^2\right) + \frac{4\nu}{h^2} K\epsilon \left( \frac{e^{\mu\tau}}    {e^{\frac{\mu}{2}\tau} - 1}     \right).
\end{eqnarray*}
 Taking the limit supremum as $t\rightarrow\infty$ we get

$$\limsup_{t\rightarrow\infty}\|w_x(t)\|_{L^2}^2 \leq \frac{2}{\mu}\left(\frac{1}{\mu\nu^2}  \epsilon^6 + \frac{2h^2}{\nu} R_2^2\, \epsilon^2\right) + \frac{4\nu}{h^2} K\epsilon \left( \frac{e^{\mu\tau}}    {e^{\frac{\mu}{2}\tau} - 1}     \right).$$  We let $\epsilon\rightarrow 0$ to obtain

$$\limsup_{t\rightarrow \infty} \|w_x(t)\|^2_{L^2} = 0, \quad \mbox{ hence,} \quad  \lim_{t\rightarrow \infty} \|w_x(t)\|^2_{L^2} = 0.$$

\subsection{The KSE with feedback control: numerical results}

The stabilization of 1D KSE has been addressed in several earlier works, for example in \cite{Ahuja2009}, \cite{Christofides2000}, and \cite{Lee_Tran_2005}, in which the common starting point is a reduced-order system  that can accurately describe the dynamics of the KSE. Then, from this resulting reduced-order system, the feedback controller can readily be designed and synthesized by taking advantage of the reduced-order techniques, called the approximate inertial manifold  and the proper orthogonal decomposition methods.  Several other works also addressed the issue of selecting the optimal actuator/sensor placement so that the desired control energy budget is achieved with minimal cost, (see for example \cite{Lou_Christofides_2002} and references therein).  

In this section we present the numerical results for the  feedback control algorithm \eqref{KSE-ih}. We illustrate through numerical simulations, the application of the proposed new method for stabilizing the unsteady zero solution of the Kuramoto-Sivashinksy equations. As a simple initial test case, we consider the example where the actuators are taken to be the first $m$ modes, where $m$ depends on the number of unstable modes and the control inputs prescribe the amplitude of the modes.  We also simulate the control algorithm using determining volume elements (local spatial averages) and as well as the determining nodal values.  

We utilize  the standard  Exponential Time Differencing fourth-order Runge-Kutta (ETDRK4 ) method and extend the algorithm to accommodate the  feedback control term.  This exponential time differencing scheme was originally derived by Cox and Matthews in \cite{CoxMat2002} and was later modified by Kassam and Trefethen in \cite{Kassam_Trefethen},  treating the problem of numerical instability in  the original scheme. This overcomes a stiff type problem via the exponential time differencing, a method which uses an idea similar to the method of integrating factor.   The implementation of ETDRK4 for the KSE equation was presented as a simple example in \cite{Kassam_Trefethen}.  We have adapted similar code for a fixed computational domain, incorporated the parameter $\nu$ as the system parameter,  and  incorporated the feedback control term appropriately.

\begin{enumerate}
\item {\it Case 1:  Globally controlled KSE with finite modes}

For a  motivating simple example, we implement the proposed feedback control system in \cite{AzTi13} using finite modes  to stabilize the steady state solution ${\bf v}=0$ of \eqref{KSE}.  We recall that the control algorithm is given as in \eqref{KSE-ih} where the
interpolant operator $I_h$ acting on $u$ is defined as in \eqref{ih-modal}.

\begin{figure}[h!]  
\centering																				
\parbox{5cm}{
\includegraphics[width=5.5cm]{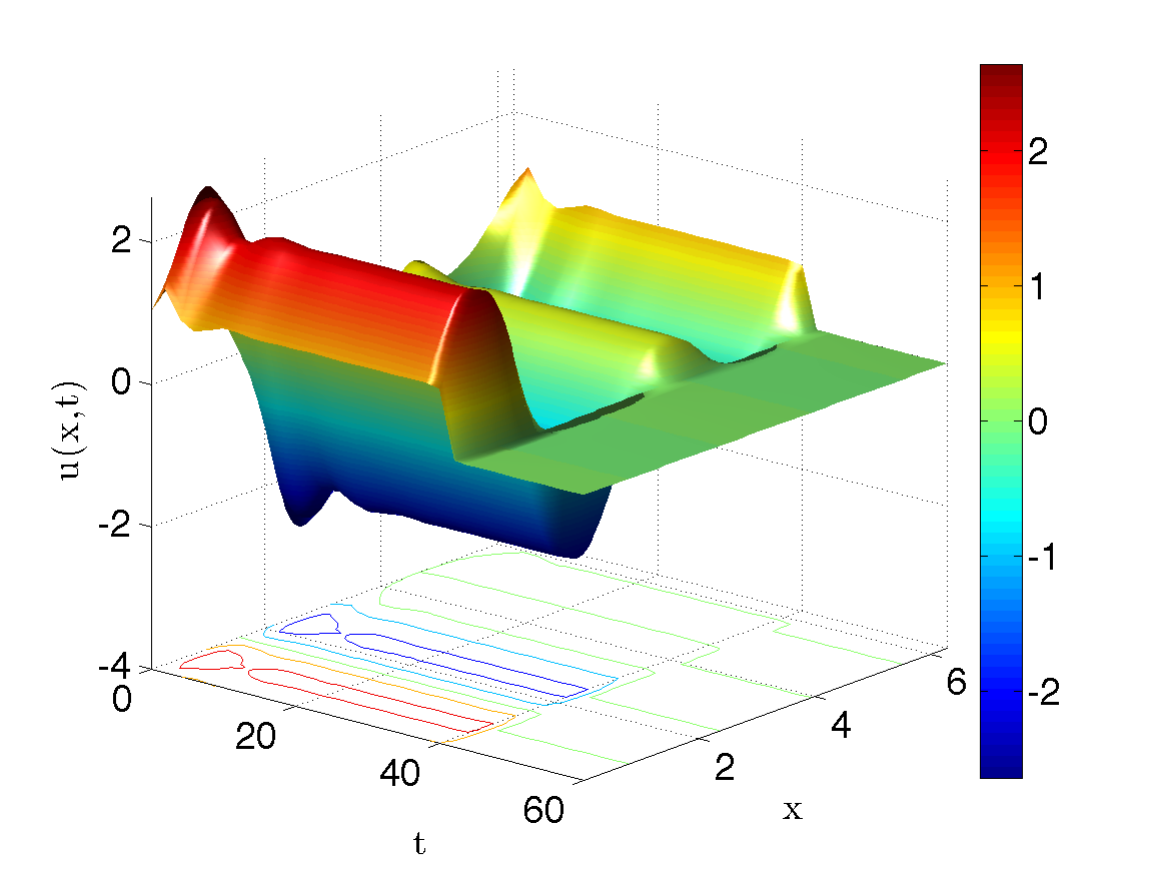}
\label{KSE1A} 																			
}
\qquad
\begin{minipage}{5cm} 																	
\includegraphics[width=5.5cm]{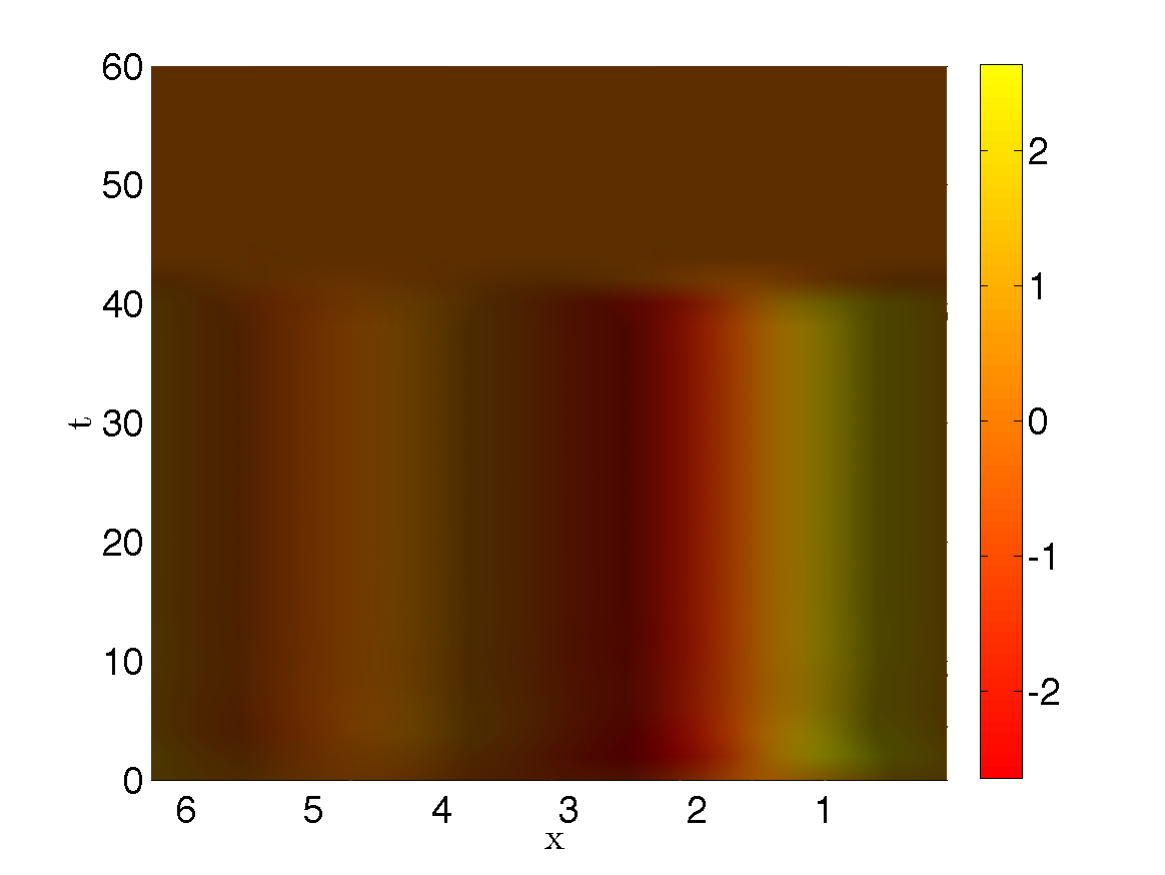}
\label{KSE1B}																			
\end{minipage}
\caption{(a) Open-loop profile showing instability of the $u(x,t)=0$ steady state solution for $0<t<40$ for $\nu = 4/15  < 1$, then the feedback control with $\mu=20$ is turned on for $t>40$ which exponentially stabilizes the system. (b)  Top view profile of $u(x, t)$.}
\end{figure}	

Here we do not assume any symmetry on the initial data except that its spatial average is zero.
Designing the feedback control based on the first $m$, $L-$periodic Fourier modes, although not practical,  is a common practice for testing grounds of the control method (for example, see \cite{Ahuja2009, Christofides2000} and references therein).  Here we present a similar initial test, that is, where the feedback control is given simply by \eqref{ih-modal}.  We are able to implement this by modifying the definition of the operator ${\bf L}$ in \cite{Kassam_Trefethen}
to accommodate the feedback control term $-\mu I_h$ and design a routine for the feedback control scheme that depends on the number of unstable modes and the type of interpolation operator $I_h$.   In Figure 4a and 4b we illustrate the solution to controlled problem using ETDRK4.  He we start with the initial condition $u_0(x) = \cos(x).*(1+\sin(x))$.  The figure in 4a illustrates the open-loop profile showing instability of the $u(x,t)=0$ steady state solution for $0<t<40$ for $\nu = 4/15  < 1$, then the feedback control with $\mu=20$ is turned on for $t>40$ which exponentially stabilizes the system. Figure 4b shows the top view profile of $u(x, t)$.

\item {\it Case 2: Globally stabilizing KSE by implementing finite volume elements}

We consider  the proposed feedback control algorithm in \eqref{KSE-ih} to stabilize the steady state solution of the 1D KSE ${\bf v}\equiv0$.
where the
interpolant operator acting on $u$, $I_h(u)$,  is defined as a slight modification of \eqref{ih-fv} as follows,

\begin{equation}
I_h(u)= \mathcal{I}_h(u) -\ang{ \mathcal{I} _h (u) },
\label{ihhh}
\end{equation}
where $I(u) = \sum_{j=1}^{N} \bar{u}_k \, \chi_{J_k \strut}(x)$, and $\displaystyle \ang{\mathcal{I}_h(u)} = \frac{1}{L}\int_0^L \mathcal{I}_hu(x,t) \, dx $.

We implement this proposed control algorithm  by modifying the main time-stepping loop via the 4th order Runge-Kutta in section 4 of \cite{Kassam_Trefethen} to accommodate the feedback control term $-\mu I_h$ and designing a function for the sensor/actuator placements.



We denote by $t_c$  the time when the feedback control is turned on.   Figure 5  illustrates the solution to controlled problem when the control is turned on at $t=0$.  He we used the  initial condition
\begin{equation}\label{e:KSE_init}
u_0(x) = \left(2.5/\sqrt{5}\right) \sum_{n=1}^5 \left(\sin (nx-n\pi) + \cos (nx-n\pi)\right).
\end{equation} The closed-loop profile shows exponential stabilization of the $u(x,t)=0$ steady state solution.  The number of controls is $NC=4$, which is proportional to the number of unstable modes.

\begin{figure}[h!]  
\centering																				
\parbox{5cm}{
\includegraphics[width=5.5cm]{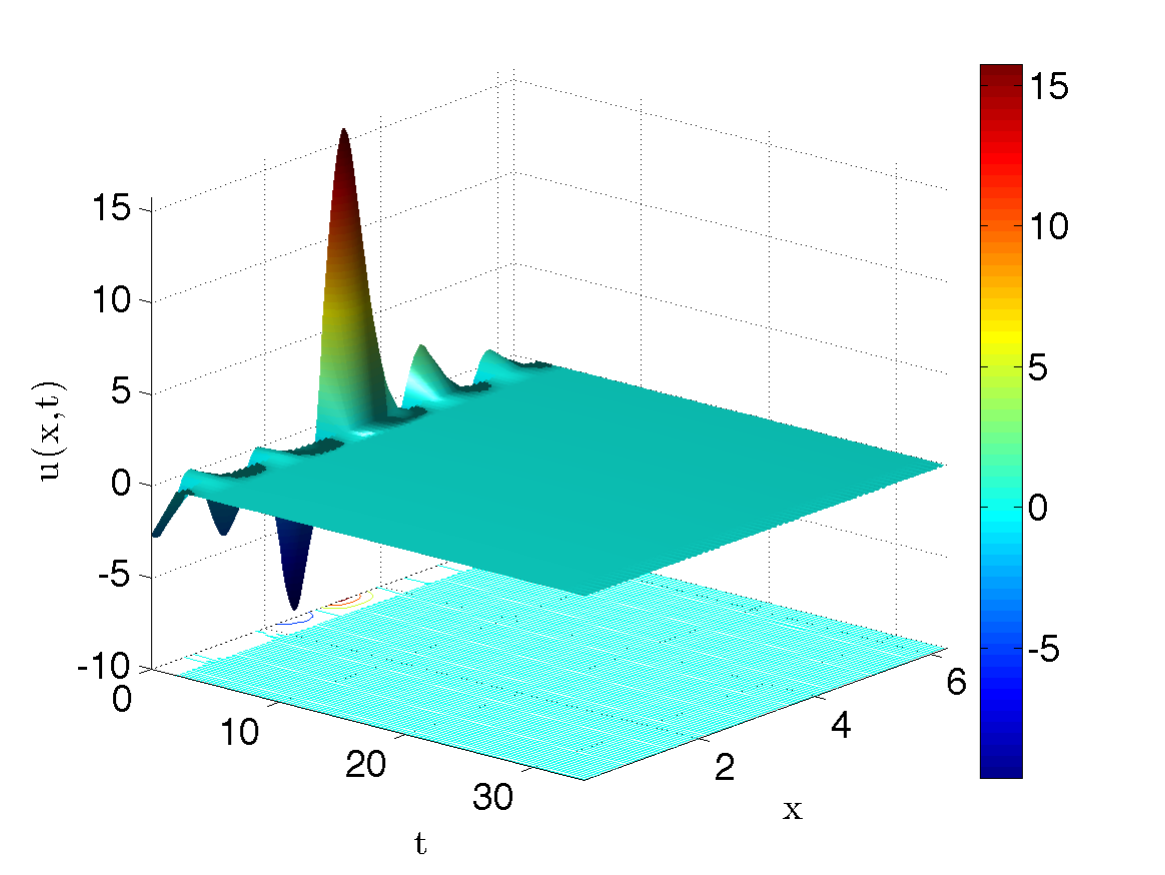}
\label{KSE2A1} 																			
}
\qquad
\begin{minipage}{5cm} 																	
\includegraphics[width=5.5cm]{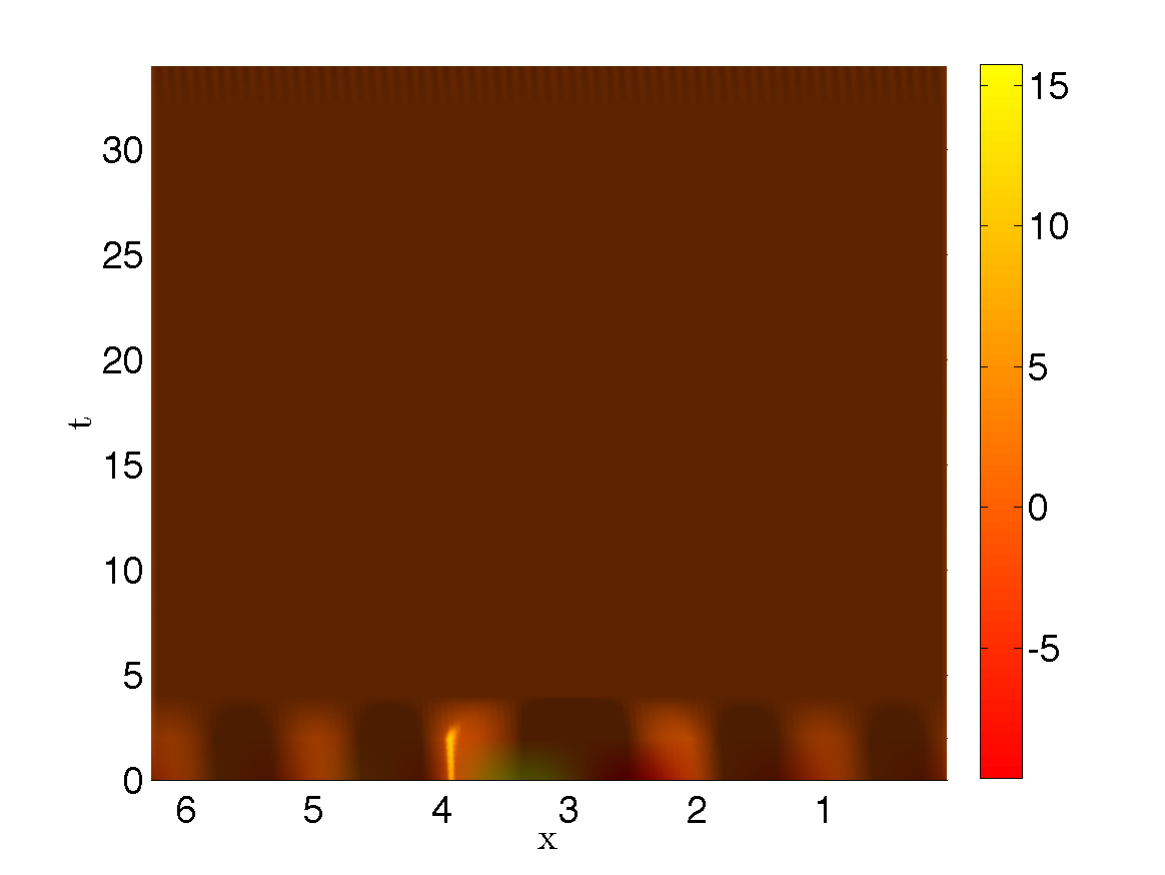}
\label{KSE2A2}																			
\end{minipage}
\caption{(a) Closed-loop profile showing  fast stabilization of the $u(x,t)=0$ steady state solution for $\nu = 4/20  < 1$, and with $\mu=20$.  (b)  Top view profile of $u(x, t)$.}
\end{figure}	

\item {\it Case 3: Globally controlled KSE with interpolant operator based on nodal values}

As we have mentioned earlier, designing the feedback control based on the first $m$, $2\pi$-periodic Fourier modes, or based on determining volume elements may not be as practical to implement in industrial setting. This is because one would require a proportional amount of sensors or controllers as the number of  grid points  used in the computer simulations.    Here we present a similar application of the control algorithm but where the feedback control uses an interpolant operator based on nodal values as defined in \eqref{ih_nodal-fv}.  The amount of physical sensors or controllers one needs is proportional to the number of unstable modes for  given parameter values.  We are able to implement this simply by modifying the feedback control routine that uses the value of the function in the middle of each subintervals instead of taking the values at every discretized spatial points and then averaging them as in the previous example with the determining volume elements.



Using similar initial condition, $u_0 = 1e^{-10}\cos x (1 + \sin x)$,  the same number of controllers $NC = 4$ and relaxation parameter $\mu = 20$ which is turned on at  $t_c=40$, we can see in Figure 6, that the film height starts to destabilize around $t = 32$ and then once  feedback control is turned on at  time $t_c=40$ it stabilizes exponentially to zero again. We recall that the simulation time is truncated to give a clear picture of the stabilization process.  We summarize our numerical experiments in Table 2.

\begin{table}
\centering
    \begin{tabular}{cclllllll}
    \toprule
    Figure & \# Actuators & $\mu$ & $\nu$ &  $t_c$ & Interpolant Operator\\
    \midrule
    2 & 0   &  0     &  1.1 & 0& \\
    3 & 0   &  0     &  4/15 & 0& \\
    4 & 4   &  20     &  4/15 &0& {\it Fourier modes} \\
    5 & 4 &  20 &  4/20  &  0& {\it finite volume}\\
    6 & 4 &  20 &  4/20  &40& {\it nodal values} \\
    \bottomrule
    \end{tabular}
    \caption{Model parameters and type of interpolant operator for the un-controlled and controlled 1D Kuramoto-Sivashinksy equations}
    \end{table}

\end{enumerate}

\begin{figure}[h!]  
\centering																				
\parbox{5cm}{
\includegraphics[width=5.5cm]{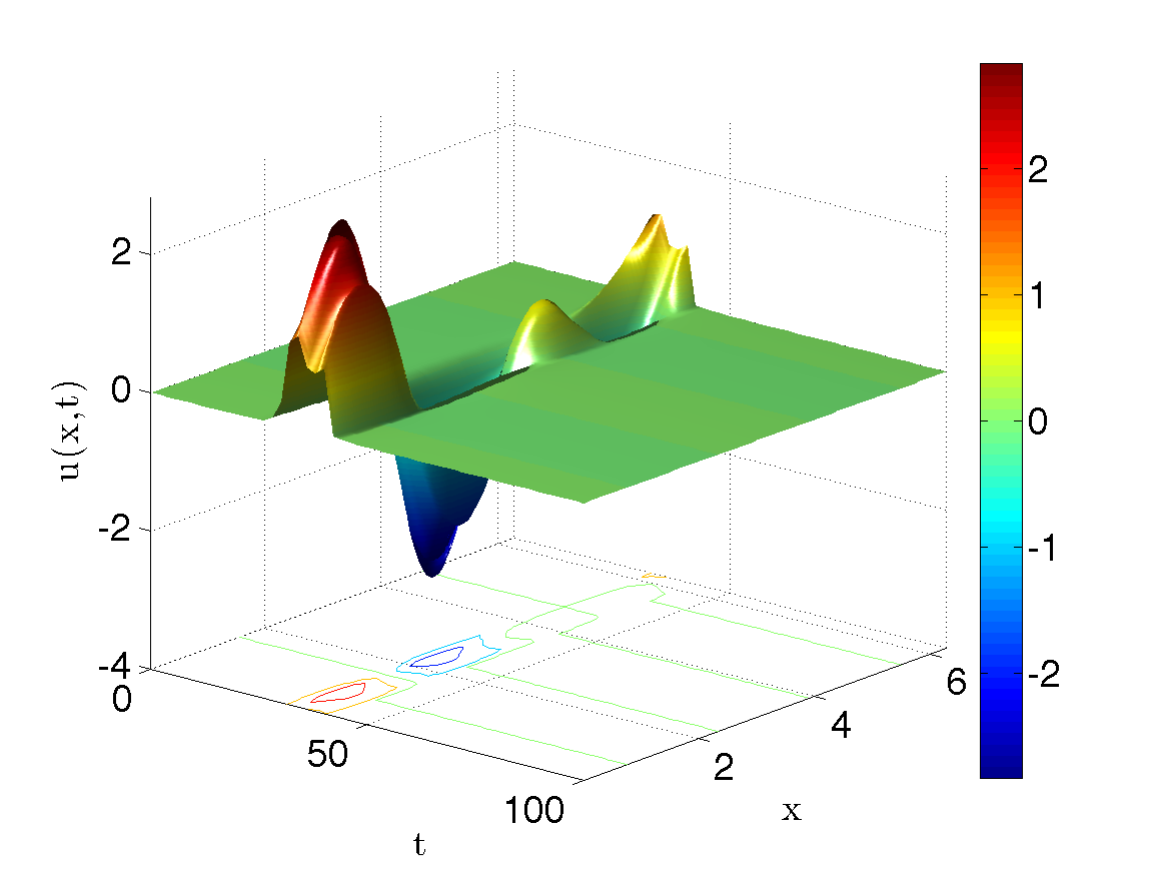}
\label{KSE2B3A} 																			
}
\qquad
\begin{minipage}{5cm} 																	
\includegraphics[width=5.5cm]{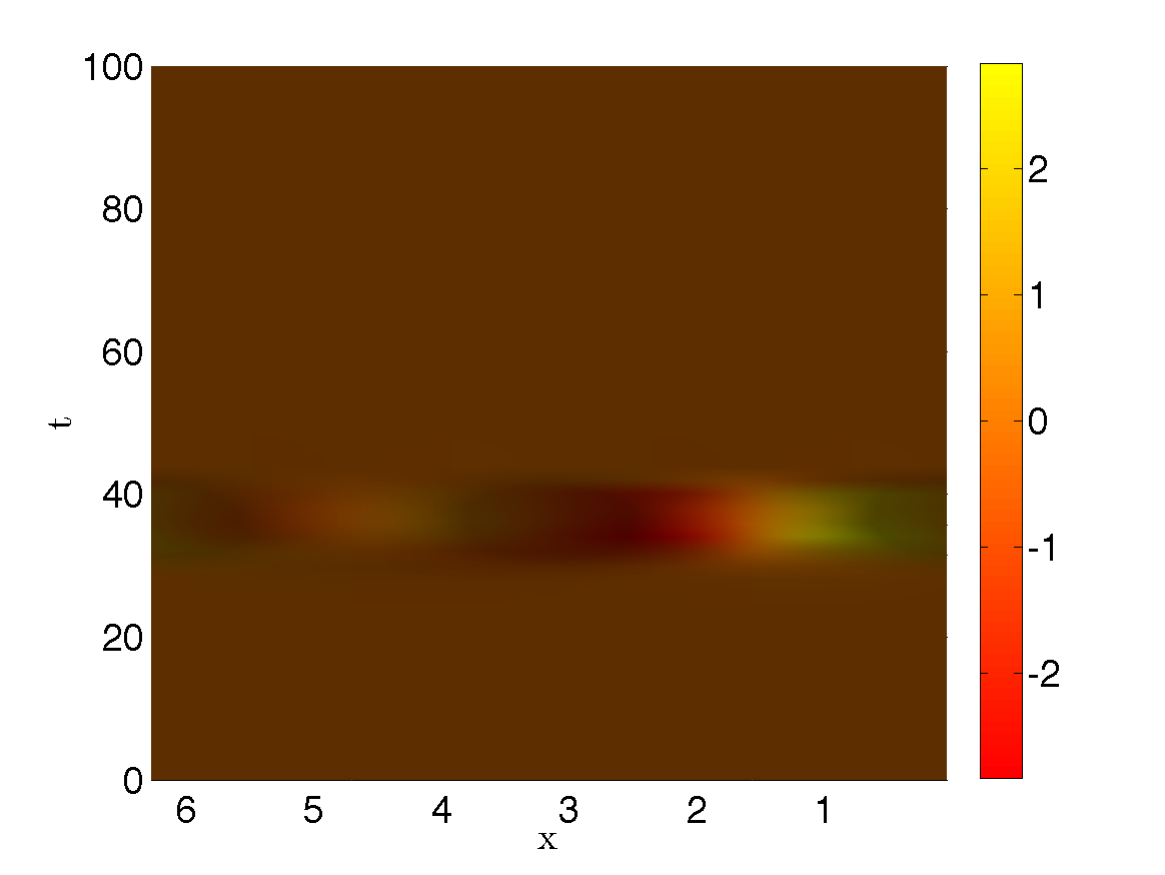}
\label{fKSE2B3B}																			
\end{minipage}
\caption{(a) With $u_0 = 1e^{-10}\cos x (1 + \sin x)$, the film height starts to destabilize around $t = 32$ and then once  feedback control is turned on at  $t_c=40$, the solution stabilizes to zero again.  (b)  A top view of the controlled profile.}
\end{figure}	

\section{Predictive control of catalytic rod with or without uncertainty variables}

We illustrate the application of the proposed method for control of nonlinear parabolic PDE system to some non-traditional feedback control case studies introduced  in \cite{Christofidesbook}.
 The setting involves  a long thin rod in a reactor where  pure species $A$ is fed into the system and a zeroth-order exothermic catalytic reaction of the form $A\rightarrow B$ takes place on the rod. Their basic assumptions in deriving the model are:
 \begin{enumerate}
 \item The reaction is assumed to be exothermic. A cooling medium is in contact with the rod to decrease the temperature of the rod.
 \item The rod has a constant density, heat capacity and constant conductivity.
 \item There is a constant temperature at both ends of the rod.
 \item  There is  unlimited supply of species A in the furnace.
 \end{enumerate} The model that describes the evolution of the dimensionless rod temperature $u(x,t)$ in the reactor  as described in \cite{Christofidesbook} is written as follows

\begin{equation}\label{cat_rod}
\frac{\partial u}{\partial t} = \frac{\partial^2 u}{\partial x^2} + \beta_T e^{-\frac{\gamma}{1 + u}} + \beta_U \left( b(x) q(t) - u  \right) - \beta_T e^{-\gamma},
\end{equation}
subject to homogeneous Dirichlet boundary conditions:

\begin{equation}
u(0,t) = 0, \quad u(\pi, t) = 0,
\end{equation}
and initial condition:
\begin{equation}
u(x,0) = u_0(x),
\end{equation}
where $\beta_T$ denotes a dimensionless heat of reaction, $\gamma$ denotes a dimensionless activation energy, $\beta_U$ denotes a dimensionless heat transfer coefficient, and $q(t)$ the manipulated input (through the cooling medium), with $b(x)$ the actuator distribution shape  function which was taken to be $b(x) = \sqrt{\frac{2}{\pi}}\sin(x)$ in \cite{Christofidesbook} chosen in order to supply maximum cooling in the middle of the rod.

\begin{enumerate}

\item {\it Case 1: Uncontrolled catalytic rod}

We start by taking the same typical values of the model parameters used in \cite{Christofidesbook}
\begin{equation}\label{cat-param}
\beta_T = 50, \quad \beta_U = 2, \quad \gamma = 4,
\end{equation}
and show that for these model parameters, the steady state solution $u(x,t) = 0$, for \eqref{cat_rod}  is unstable.   Starting with an initial data with small perturbation near zero, the temperature evolves to another stable steady state  where the temperature profile has a hot-spot in the middle.    We run the simulation with the initial condition $u_0(x) = 1e^{-3}\sin(2x)$ on the spatial interval $[0 , \pi]$ and time interval $[0, 6]$.  We obtain the following results illustrated in Figure 7.  The axes are in units of $\Delta x = \pi/20$ and $\Delta t = 6/1000$.

\begin{figure}[h!]  
\centering																				
\parbox{5cm}{
\includegraphics[width=5.5cm]{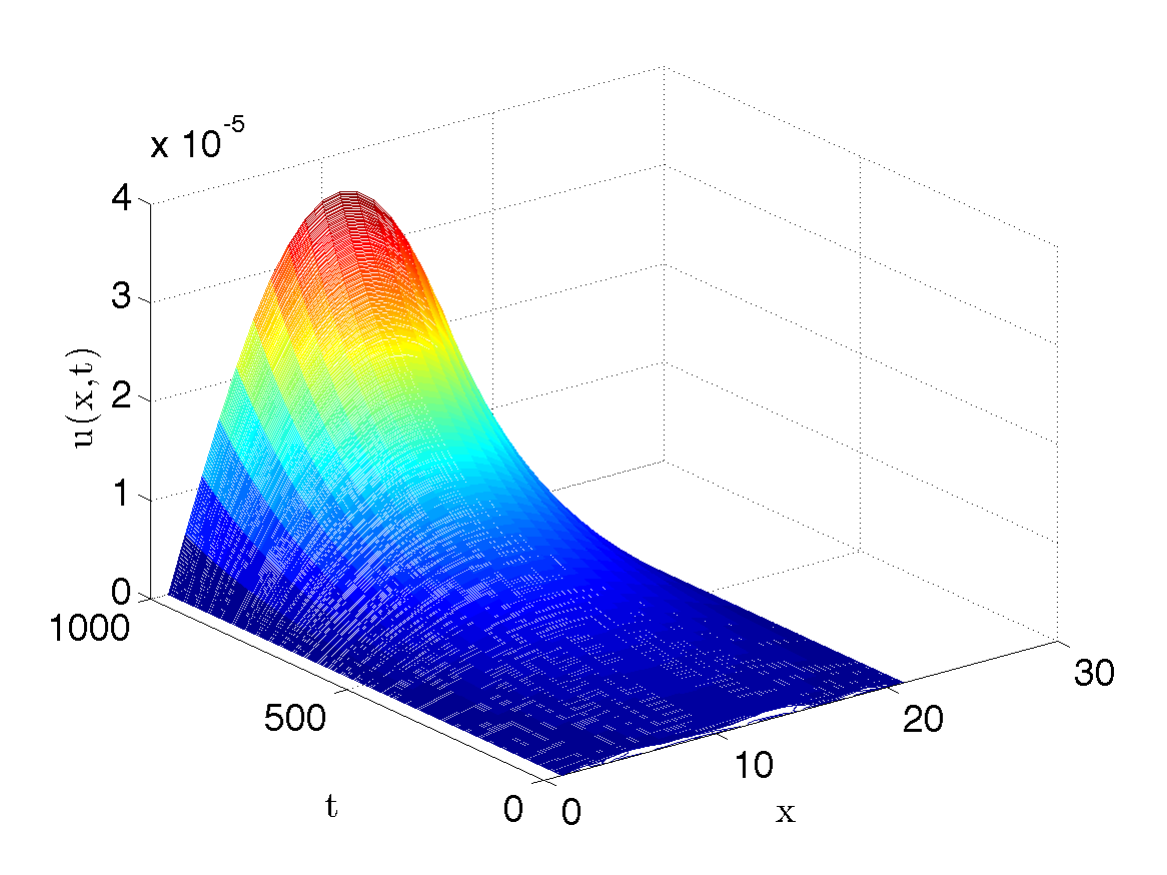}
\label{CAT1A} 																			
}
\qquad
\begin{minipage}{5cm} 																	
\includegraphics[width=5.5cm]{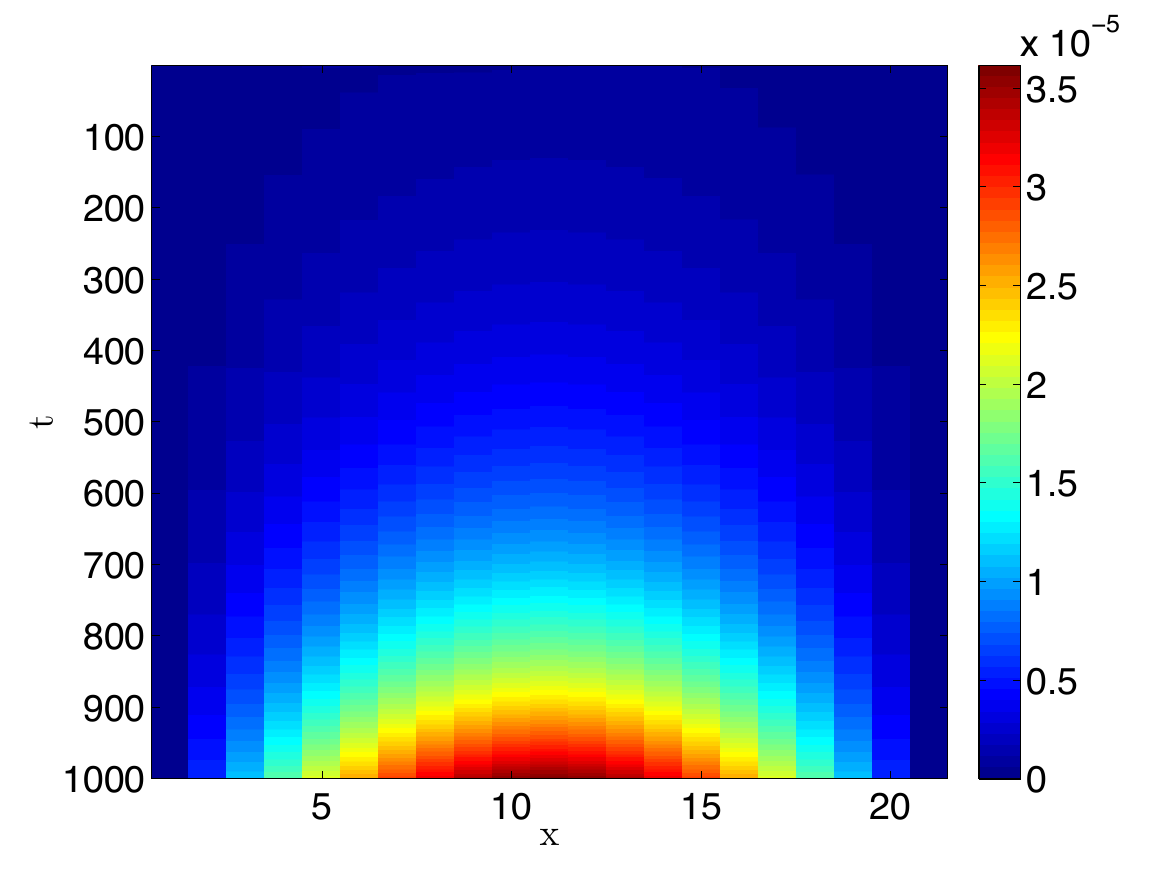}
\label{CAT1B}																			
\end{minipage}
\caption{(a) Open-loop profile showing instability of the $u(x,t)=0$ steady state solution. (b)  Top-view of  $u(x,t)$.}
\end{figure}	

\item{ \it Case 2:  Globally controlled catalytic rod}

We now apply the feedback control scheme proposed earlier for the catalytic rod problem which takes the following form:

\begin{equation}\label{cat_rod_ih}
\frac{\partial u}{\partial t} = \frac{\partial^2 u}{\partial x^2} + \beta_T e^{-\frac{\gamma}{1 + u}} + \beta_U \left( -\mu I_h(u) - u  \right)   - \beta_T e^{-\gamma}.
\end{equation}
\noindent
For the given parameters for the catalytic rod problem in \eqref{cat-param}, we observe one unstable mode and so we supply our control algorithm with the number of controller $NC=1$.  We put one actuator in the middle of the rod at $x = \pi/2$.  The interpolant operator is defined as $I_h(u) =  \bar{u} \, \chi_{[0,\pi] \strut}(x)$, where $\bar{u}$ is the spatial average of $u(x,t)$ on the interval $[0, \pi]$.  Under this feedback control scheme, using the same initial condition,  we observed global stabilization of the trivial steady state solution as shown in Figure 8.

\begin{figure}[h!]  
\centering																				
\parbox{5cm}{
\includegraphics[width=5.5cm]{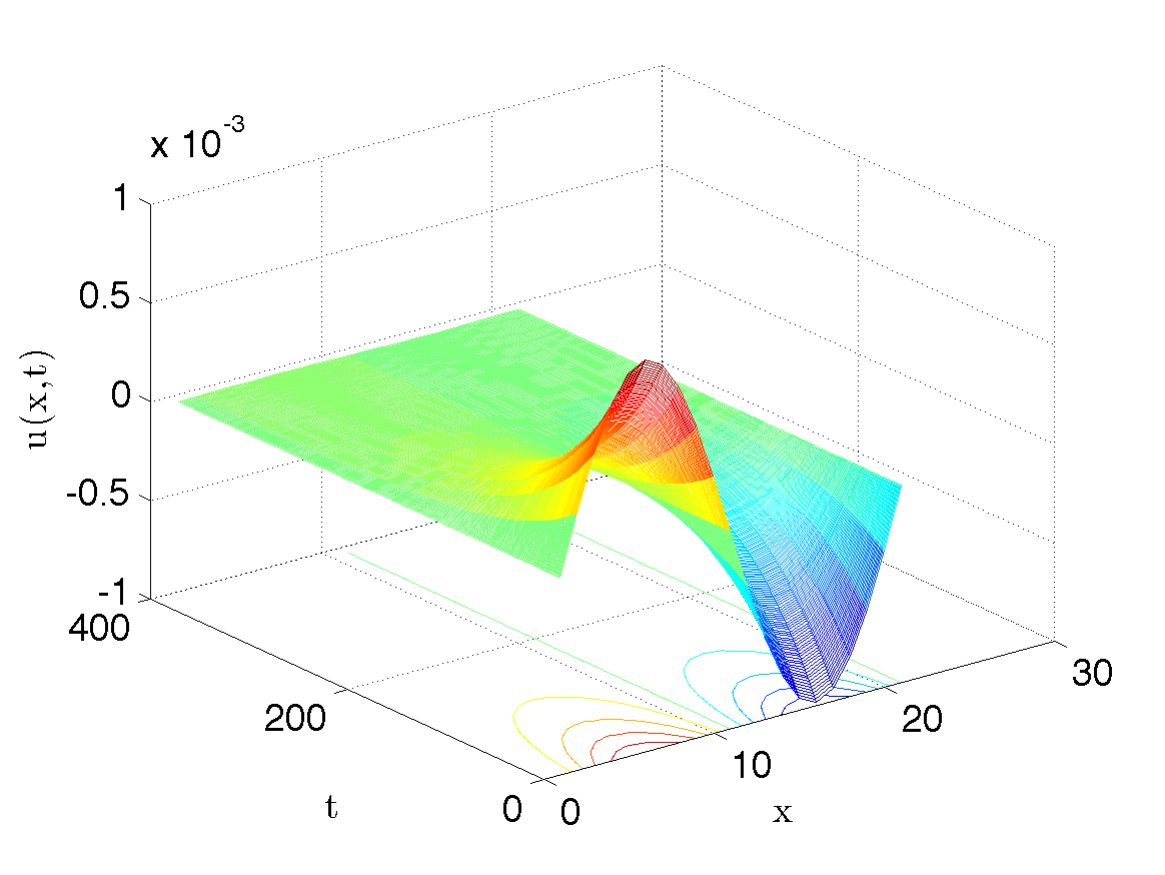}
\label{CAT2A1} 																			
}
\qquad
\begin{minipage}{5cm} 																	
\includegraphics[width=5.5cm]{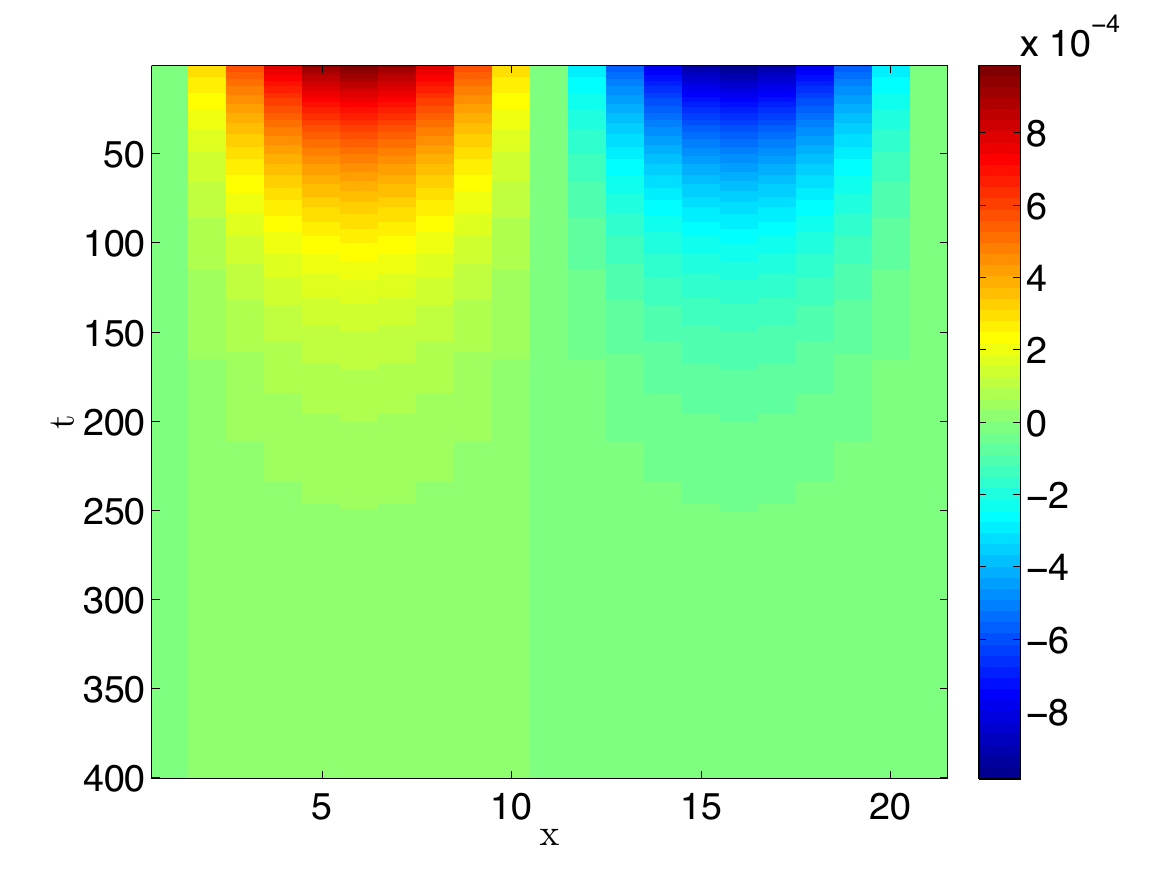}
\label{CAT2A2} 																		
\end{minipage}
\caption{(a) Closed-loop profile showing stabilization to $u(x,t)=0$ steady state solution. (b)  Top-view.}
\end{figure}	

\item {\it Case 3: Application to catalytic rod with uncertainty}

Motivated by an example in \cite{Christofidesbook}, we investigate the performance of the feedback control algorithm when there is a time varying uncertainty in some of the parameters in the model equation.  Our numerical simulation shows that our feedback control algorithm is able to regulate the temperature profile in the rod through the manipulation of the temperature of the cooling medium even in the presence of time-varying uncertainty in the heat of the reaction $\beta_T$.  The author in \cite{Christofidesbook} proposed a procedure for the synthesis of the robust controllers that achieve arbitrary degree of asymptotic attenuation of the effect of the parameters with uncertainty on the output based on the construction of higher dimensional approximation of the state slow-variables subsystems stemming from the concept of inertial manifold.   Here we apply a similar study using the linear feedback control algorithm.  
For simplicity, following \cite{Christofidesbook}, we used $\beta_T = \overline{\beta}_{T} + \theta(t)$, where $\overline{\beta}_{T} = 50$ and $\theta(t) = \sin(0.524 t)$.  The location of the actuator is at $x = \pi/2$.  In the presence of this kind of uncertainty, our control algorithm is able to regulate the temperature with an error related to the size of the uncertainty in the model parameters.   Our results are shown in Figure 9 for the case where the initial condition has amplitude of $A = 1e^{-10}$.  We observe the eventual stabilization to zero.

\begin{figure}[h!]  
\centering																				
\parbox{5cm}{
\includegraphics[width=5.5cm]{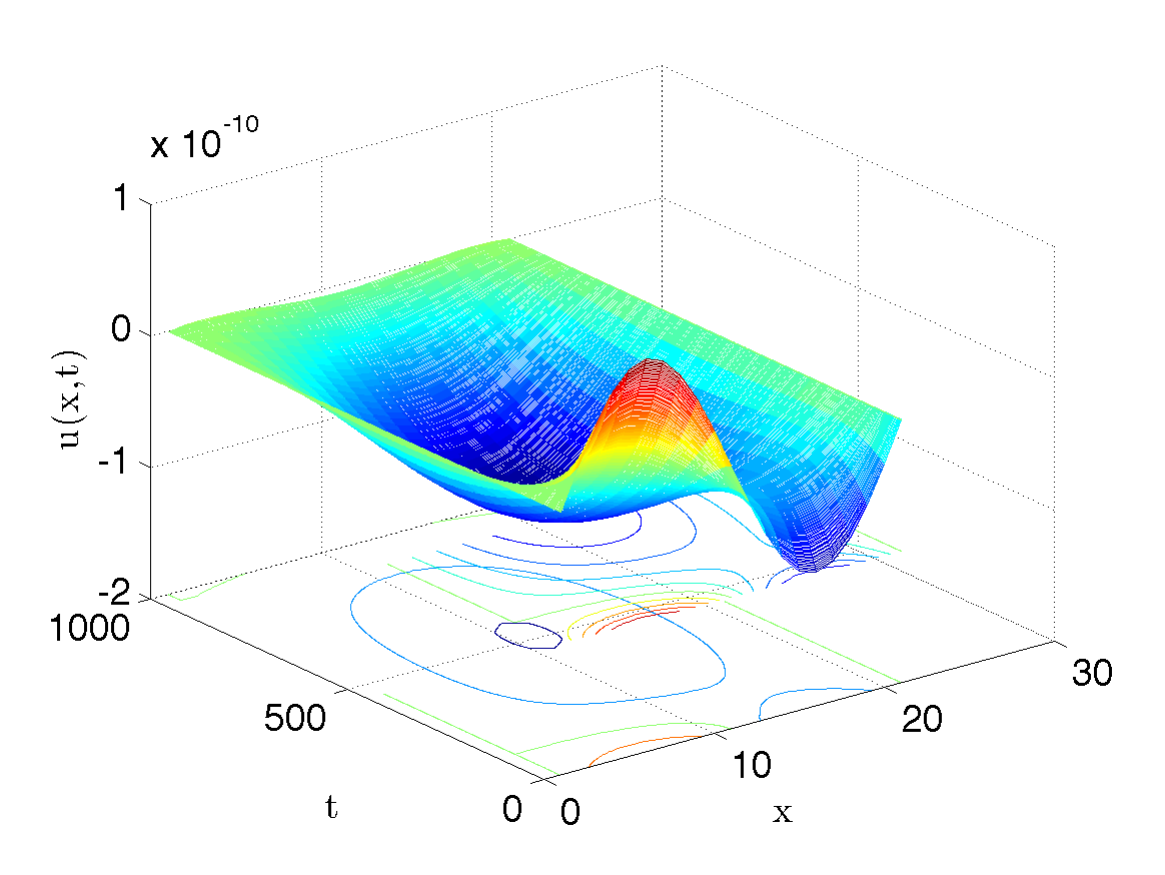}
\label{CAT3B1} 																		
}
\qquad
\begin{minipage}{5cm} 																	
\includegraphics[width=5.5cm]{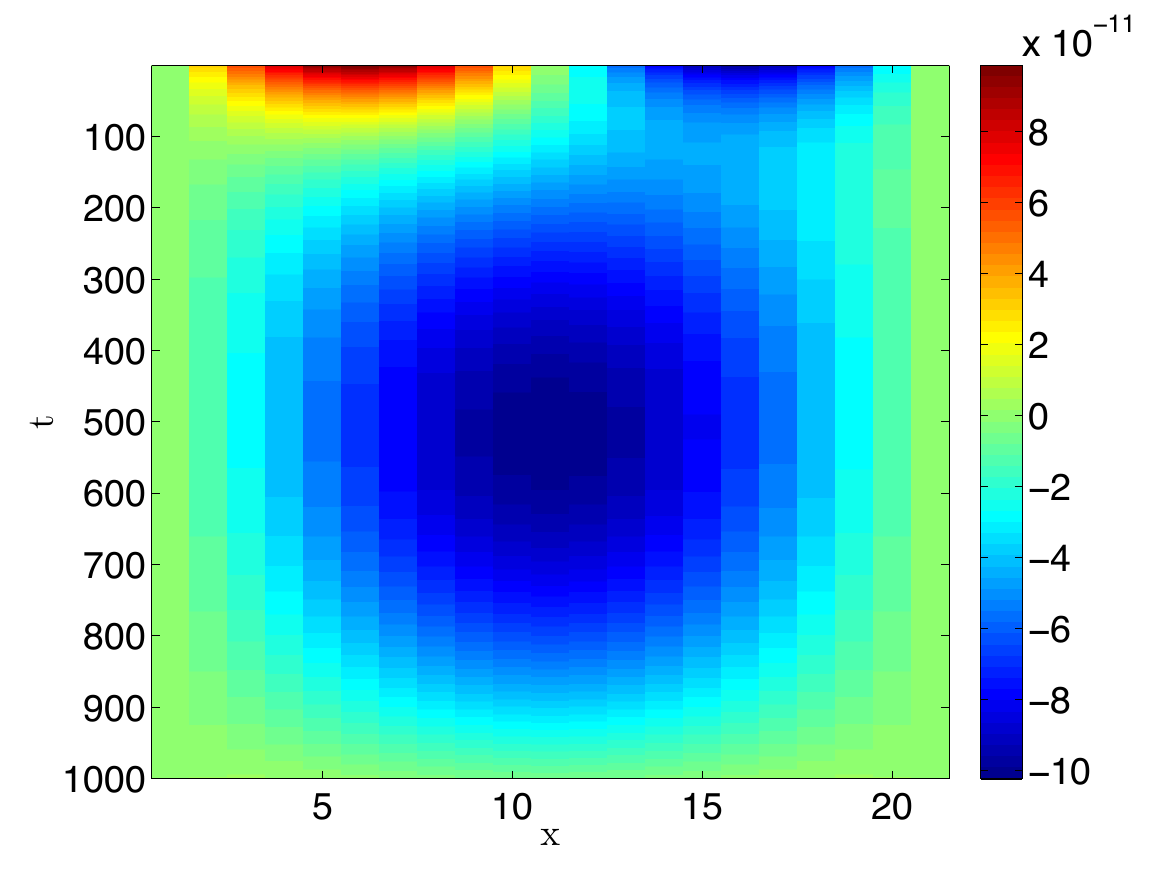}
\label{CAT3AB2} 																			
\end{minipage}
\caption{(a) Closed-loop profile showing eventual stability. (b)  Top-view.}
\end{figure}	


\item {\it Case 4: Nodal-point observational measurements}

For our last numerical test, we repeat Case 2 but with the feedback control scheme is given by
\begin{equation}\label{cat_rod_ih}
\frac{\partial u}{\partial t} = \frac{\partial^2 u}{\partial x^2} + \beta_T e^{-\frac{\gamma}{1 + u}} + \beta_U \left( -\mu I_h(u) - u  \right)   - \beta_T e^{-\gamma},
\end{equation}
\noindent
with $I_h(u) =  u(\pi/2)\chi_{[0,\pi]}$.  So the actuator and sensor are both located in the middle of the rod at $x = \pi/2$.   We have observed similar behavior as in the case of the finite volume case.   For this reason we do not present the figures.

We summarize our numerical experiments in Table 3.

\begin{table}[h!]
\centering
    \begin{tabular}{cclllllllllll}
    \toprule
 Figure & \# Actuators & $\mu$ & $\nu$ & $\beta_T$ &$\beta_U$ &$\gamma$ & interpolant operator\\
    \midrule
    7 & 0   &  0     &  1 &   50&2.0 &4.0&  \\
    8 & 1 &  30 &  1 &   50 &  2.0&4.0&finite volume\\
      9 & 1 &  30 &  1 &varying  &  2.0&4.0&finite volume\\
      similar to Fig 8 & 1 & 30 & 1& 50& 2.0&4.0& nodal values\\
    \bottomrule
    \end{tabular}
    \caption{Model parameters and type of interpolant operator for the un-controlled and controlled catalytic rod problem}
    \end{table}
\end{enumerate}
 \subsection*{Acknowledgements}
 The work of E.L.~is supported  by the ONR grant N001614WX30023.
The work of E.S.T. was supported in part by the ONR grant N00014-15-1-2333
and by the NSF grants DMS--1109640 and DMS--1109645.
 E.S.T. is also thankful to the kind hospitality  of the Department of Mathematics, United States Naval Academy where this work was initiated.

\end{document}